\documentclass[10pt]{article}
\usepackage{grffile}
\usepackage[normalem]{ulem}
\usepackage[margin=1in]{geometry}
\usepackage{amssymb, cancel, amsmath, amstext}
\usepackage[dvips]{epsfig}
\usepackage[small]{caption}
\usepackage{graphicx}
\usepackage[all]{xy}
\usepackage{mathtools}
\usepackage{hyperref}

\usepackage{caption}
\usepackage{subcaption}
\usepackage{colortbl}
\usepackage{mathrsfs}
\usepackage{tikz}
\usetikzlibrary{positioning}

\usepackage{ntheorem}

\usepackage[shortlabels]{enumitem}
\usepackage{wrapfig}
\usepackage{etex}

\usepackage{marginnote}
\usepackage{geometry}
\newtheorem{Lemma}{Lemma}[section]
\newtheorem{Theorem}[Lemma]{Theorem}
\newtheorem{Proposition}[Lemma]{Proposition}
\newtheorem{Corollary}[Lemma]{Corollary}
\newtheorem{Remark}[Lemma]{Remark}
\newtheorem{Definition}[Lemma]{Definition}

\newtheorem*{Assumption}[Lemma]{Assumption}

\newenvironment{Proof}[1][\unskip]%
 {\begin{trivlist} \item[]{\bf Proof #1. }}%
 {\hspace*{\fill}$\rule{.4\baselineskip}{.4\baselineskip}$\end{trivlist}}

{\begin{trivlist}\item[]\textbf{Acknowledgments.}}{\end{trivlist}}

\makeatletter\@addtoreset{figure}{section}\makeatother

\makeatletter \@addtoreset{equation}{section} \makeatother

\newcommand{\tx}{\tilde{x}}
\newcommand{\ty}{\tilde{y}}

\newcommand{\R}{\mathbb{R}}

\newcommand{\Z}{\mathbb{Z}}

\newcommand{\rmd}{\mathrm{d}}
\newcommand{\rme}{\mathrm{e}}
\newcommand{\rmi}{\mathrm{i}}

\renewcommand{\ker}{\mathrm{Ker}\,}

\renewcommand{\leq}{\leqslant}
\renewcommand{\geq}{\geqslant}

\def\XXint#1#2#3{{\setbox0=\hbox{$#1{#2#3}{\int}$}
\vcenter{\hbox{$#2#3$}}\kern-.5\wd0}}

\begin{document}

\title{Contact angle selection for interfaces in growing domains}

\author{ Rafael Monteiro  and  Arnd Scheel  \\[2ex]
\textit{\footnotesize University of Minnesota, School of Mathematics,  
206 Church St. S.E., Minneapolis, MN 55455, USA}} 

\date{\small \today} 

\maketitle 
\begin{abstract}
We study interfaces in an Allen-Cahn equation, separating two metastable states. Our focus is on a directional quenching scenario, where a parameter renders the system bistable in a half plane and monostable in its complement, with the region of bistability expanding at a fixed speed. We show that the growth mechanism selects a contact angle between the boundary of the region of bistability and the interface separating the two metastable states. Technically, we focus on a perturbative setting in a vicinity of a symmetric situation with perpendicular contact. The main difficulty stems from the lack of Fredholm properties for the linearization in translation invariant norms. We overcome those difficulties establishing Fredholm properties in weighted spaces and farfield-core decompositions to compensate for negative Fredholm indices. 

\vspace{\baselineskip}

\end{abstract}

\section{Introduction}
We are interested in the Allen-Cahn equation with directional quenching,
\begin{equation}\label{Allen-Cahn}
u_t  = \Delta u + \mu(x-c_x t) u  - u^3 + \alpha g(x-c_x t,u) ,\qquad (x,y)\in\R^2.
\end{equation}
Here, $u$ denotes an order parameter, and $\alpha$ is a perturbation parameter. We assume for simplicity that $\mu(x)$ and $g(x,u)$ are constant in $x>0$ (right half plane) and $x<0$ (left half plane),
\[
\mu(x)=\left\{\begin{array}{ll}+1,&x<0\\-1,&x>0\end{array}\right. ,\qquad 
g(x,u)=\left\{\begin{array}{ll}g_\mathrm{l}(u),&x<0\\g_\mathrm{r}(u),&x>0\end{array}\right.
.\]
In particular, for $\alpha=0$, the origin $u=0$ is the unique stable equilibrium in the region $x>c_xt$ and $u=\pm 1$ are the two stable equilibria in $x<c_xt$. Throughout, when referring to solutions of \eqref{Allen-Cahn}, we refer to weak solutions, understanding that solutions are smooth away from the jump $x=c_xt$ and possess continuous derivatives across the jump. 

\begin{figure}[htb]
\begin{center}
\hspace*{0.2in}\includegraphics[height=1.in]{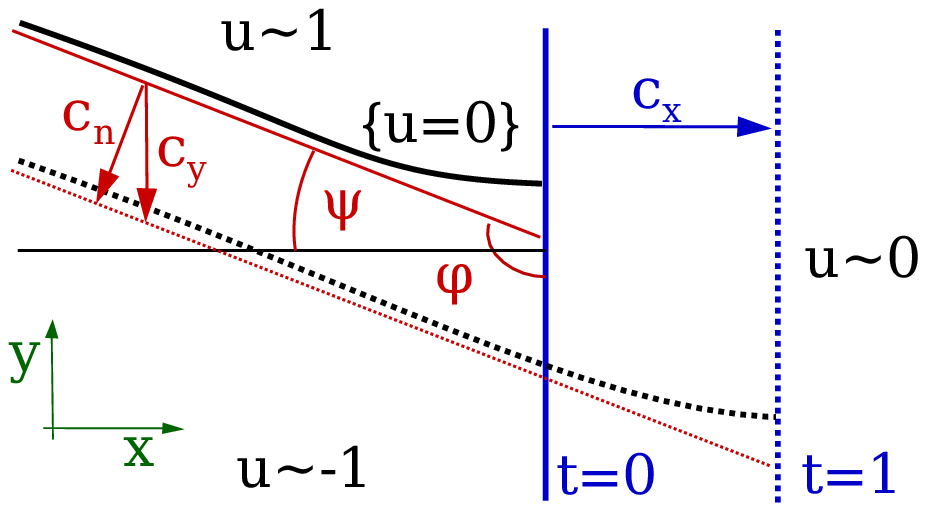} \hfill \includegraphics[height=1.in]{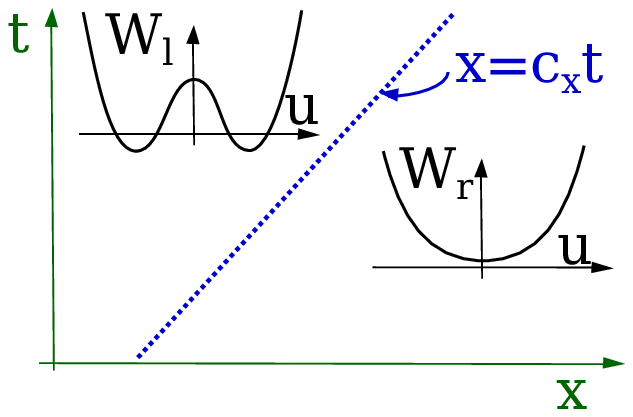} \hfill
\includegraphics[height=1.in]{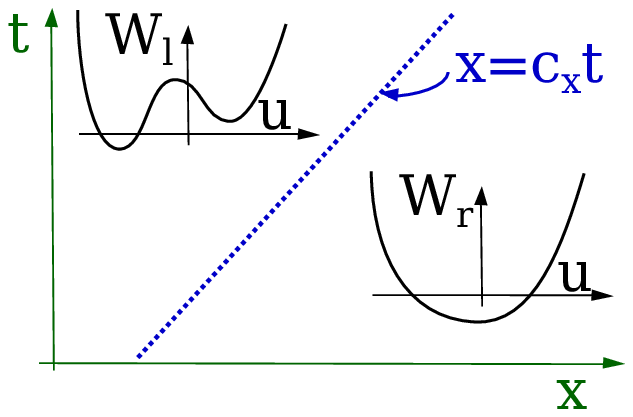}\hspace*{0.2in}
\end{center}
\caption{Schematic figure of quenching line (blue), asymptotes (red), and nodal line  of interface (black), at time $t=0$ and time $t=1$, including contact angles and interfacial speeds (left). Local potential enegies  $W$ in space-time, for $\alpha=0$, balanced, in the center figure, and for $\alpha\neq 0$, unbalanced,  on the right-hand side.\label{f:1}}
\end{figure}

We view this equation as a simple model for phase separation through directional quenching. The quenching line $\{y\in\R,x=c_xt\}$ separates bistable and monostable regions. The pointwise energy $W_\mathrm{l/r}(u)$, with 
\[
-W'_\mathrm{l}(u)=u-u^3+\alpha g_\mathrm{l}(u),\qquad 
-W'_\mathrm{r}(u)=-u-u^3+\alpha g_\mathrm{r}(u),
\]
is of double-well type in the left half-plane and convex in the right half-plane; see Figure \ref{f:1}.

The speed $c_x>0$ measures the rate of growth of the bistable region. The parameter $\alpha$ encodes possible asymmetries, that is, situations where for instance $W_\mathrm{l}(+1)\neq W_\mathrm{l}(-1)$, or $W'_\mathrm{r}(0)\neq 0$. 

Loosely speaking, we are interested in solutions where 
\begin{equation}\label{e:bc0}
u(x,y)\longrightarrow \left\{\begin{array}{ll} 
0,& x\to +\infty,\\
z_+>0,& x\to -\infty,\ y\to +\infty\\
z_-<0,& x\to -\infty,\ y\to -\infty
\end{array}\right.,
\end{equation}
where $z_\pm$ denote the local minimizers in $x<c_xt$,  $W_\mathrm{l}'(z_\pm)=0$. 

For $x$ negative, large, we then expect an interface between the regions where $u\sim z_+$ and where $u\sim z_-$, marked for instance by the level set  $\{u=0\}$. This interface typically propagates with a distinct normal speed $c_\mathrm{n}(\alpha)$, small for $\alpha\sim 0$, while the region  $x<ct$, in which the interface actually describes the dynamics, expands to the right. 

For $\alpha=0$, we previously showed that there exists a solution $u(x,y)$ that is odd in $y$ \cite{Monteiro_Scheel}. In particular, the interface mentioned above consists of the $x$-axis, $y=0$. We show here that this solution can be continued as a traveling wave, that is, a stationary solution in a frame $\tilde{x}=x-c_x t$, $\tilde{y}=y-c_y t$,  with a selected vertical speed $c_y(\alpha)$ and a selected asymptotic angle $\varphi(\alpha)$, that is, $u=0$ on $\tilde{y}=\cot(\varphi)\tilde{x}$ as $x\to -\infty$; see Figure \ref{f:1}.

Questions of directional quenching are part of a larger set of questions of self-organization in growing domains. The question of interest is how the growth process, here the movement of the quenching line $x=ct$, acts as a selection mechanism for structure, such as interfaces or more complicated patterns, in its wake. Such selection mechanisms may be exploited, in biology or engineering; see \cite{Monteiro_Scheel} for a somewhat broader overview. 

Transforming to the comoving frame, and dropping tildes for the new variables, we find the elliptic traveling-wave equation, denoting partial derivatives by subscripts, 
\begin{equation}\label{e:ell}
\Delta u +c_x u_x + c_y u_y +  \mu(x) u  - u^3 + \alpha g(x,u)=0 ,\qquad (x,y)\in\R^2.
\end{equation}
As a next step, we will make the rough formulation of boundary conditions \eqref{e:bc0} more precise. Denote by $z_\pm(\alpha)$ the unique, smooth family of zeros of $u-u^3+\alpha g_\mathrm{l}(u)$ with $z_\pm(0)=\pm1$, and by $z_0(\alpha)$ the unique zero of $-u-u^3+\alpha g_\mathrm{r}(u)$, $|\alpha|$ sufficiently small.

\begin{Definition}[Contact angle]\label{d:contact}
We say \eqref{e:ell} possesses a solution $u$ with contact angle $\varphi_*$ if $u$ possesses the limits 
\begin{equation}\label{e:bc}
\lim_{x\to +\infty} u(x,y)=z_0(\alpha), \qquad \lim_{x\to -\infty} u(x,\cot(\varphi)x)=\left\{\begin{array}{ll}
z_+(\alpha),& \varphi>\varphi_*\\
z_-(\alpha),& \varphi<\varphi_*
\end{array}\right. ,
\end{equation}
for all $0<\varphi<\pi$. We also use the deviation from a right angle, $\psi=\varphi-\frac{\pi}{2}$; see Figure \ref{f:1}.
\end{Definition}
Our main result provides the existence of solutions with a selected contact angle. 
\begin{Theorem}[Existence and contact angle selection]\label{t:1}
For $0<c_x\ll1$ sufficiently small, there exists $\alpha_0(c_x)>0$ such that for all $|\alpha|<\alpha_0(c_x)$ there exist a speed $c_y(\alpha)$ and a solution  $u(x,y;\alpha)$ to \eqref{e:ell} with contact angle $\varphi(\alpha)$. Moreover, $c_y(\alpha)$ and $\varphi(\alpha)$ are smooth with $\varphi(0)=\pi/2$, $c_y(0)=0$, and $u(x,y;\alpha)$ is smooth in $\alpha$  in a locally uniform topology, that is, considering the restriction $u|_{B_R(0)}$ to an arbitrary large ball.
\end{Theorem}

\begin{figure}
\includegraphics[width=0.16\textwidth]{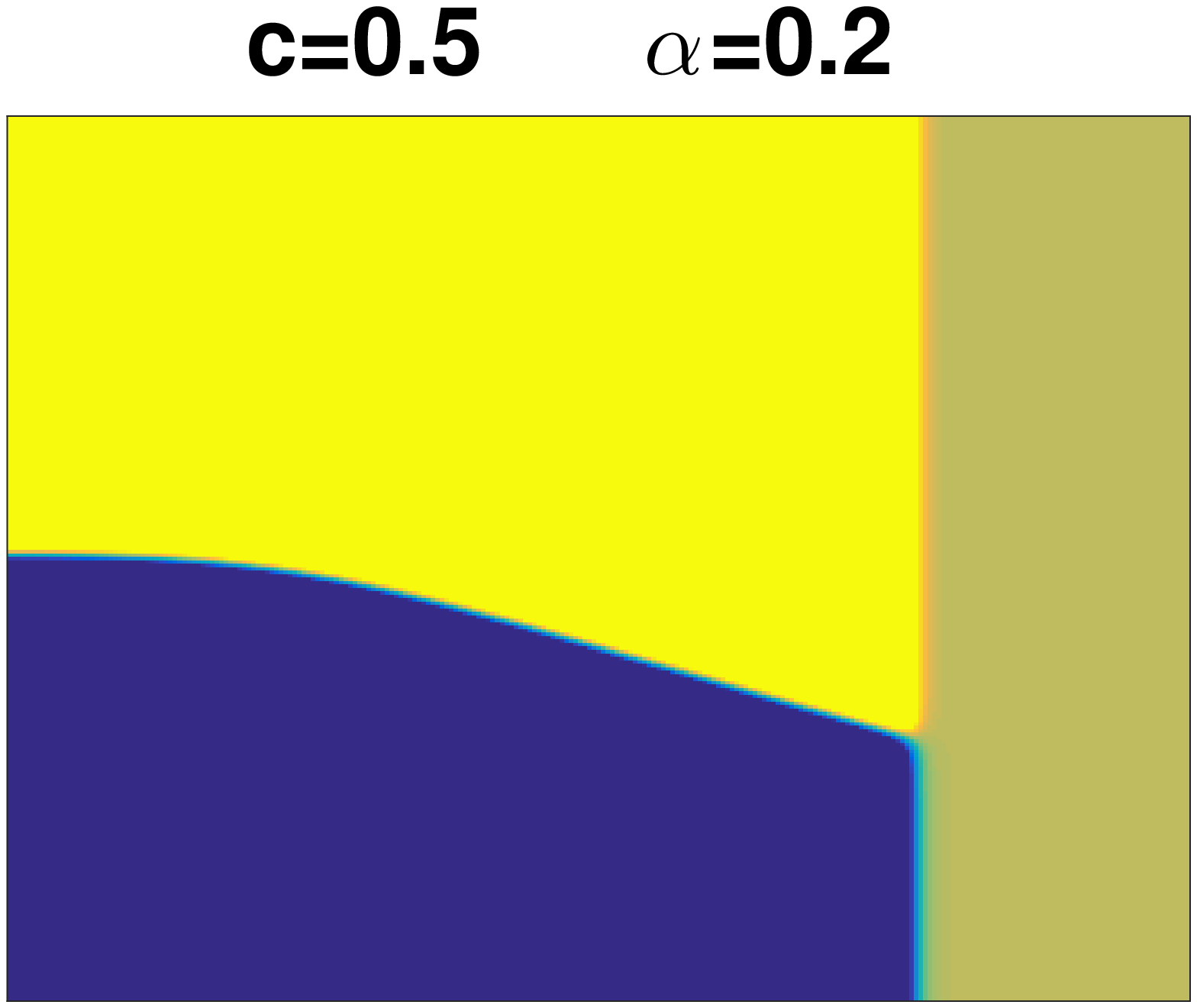}\hfill
\includegraphics[width=0.16\textwidth]{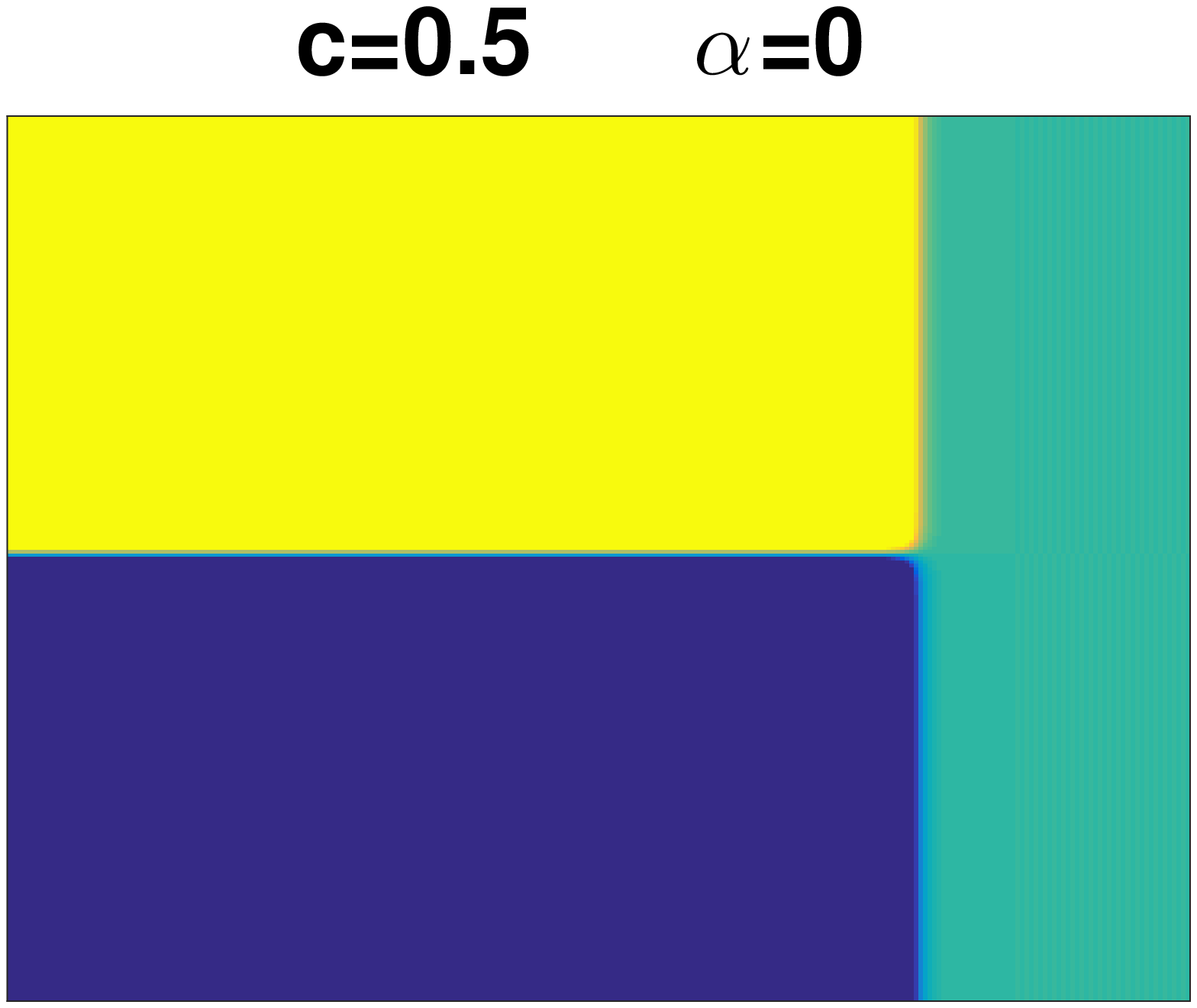}\hfill
\includegraphics[width=0.16\textwidth]{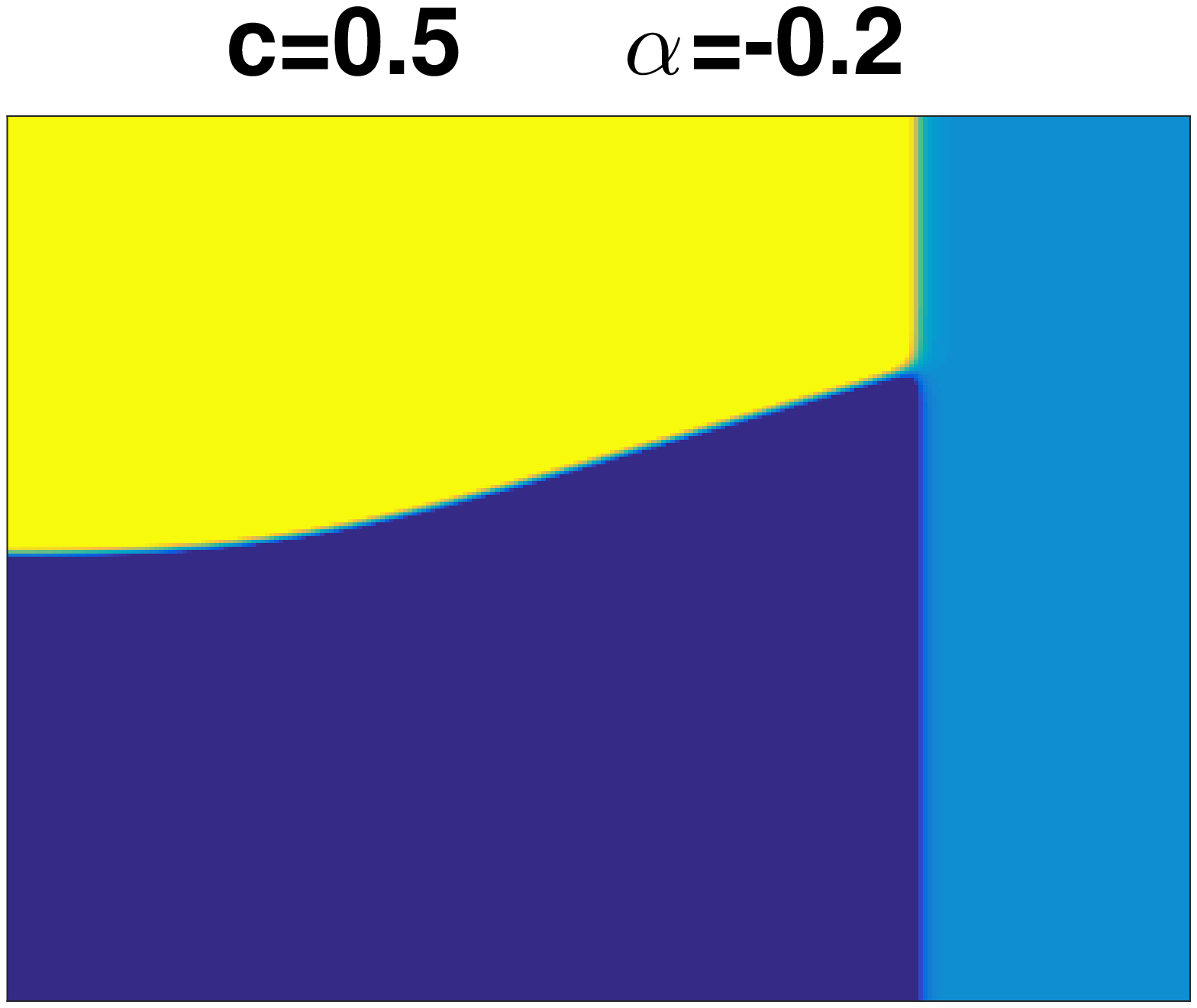}\hfill
\includegraphics[width=0.16\textwidth]{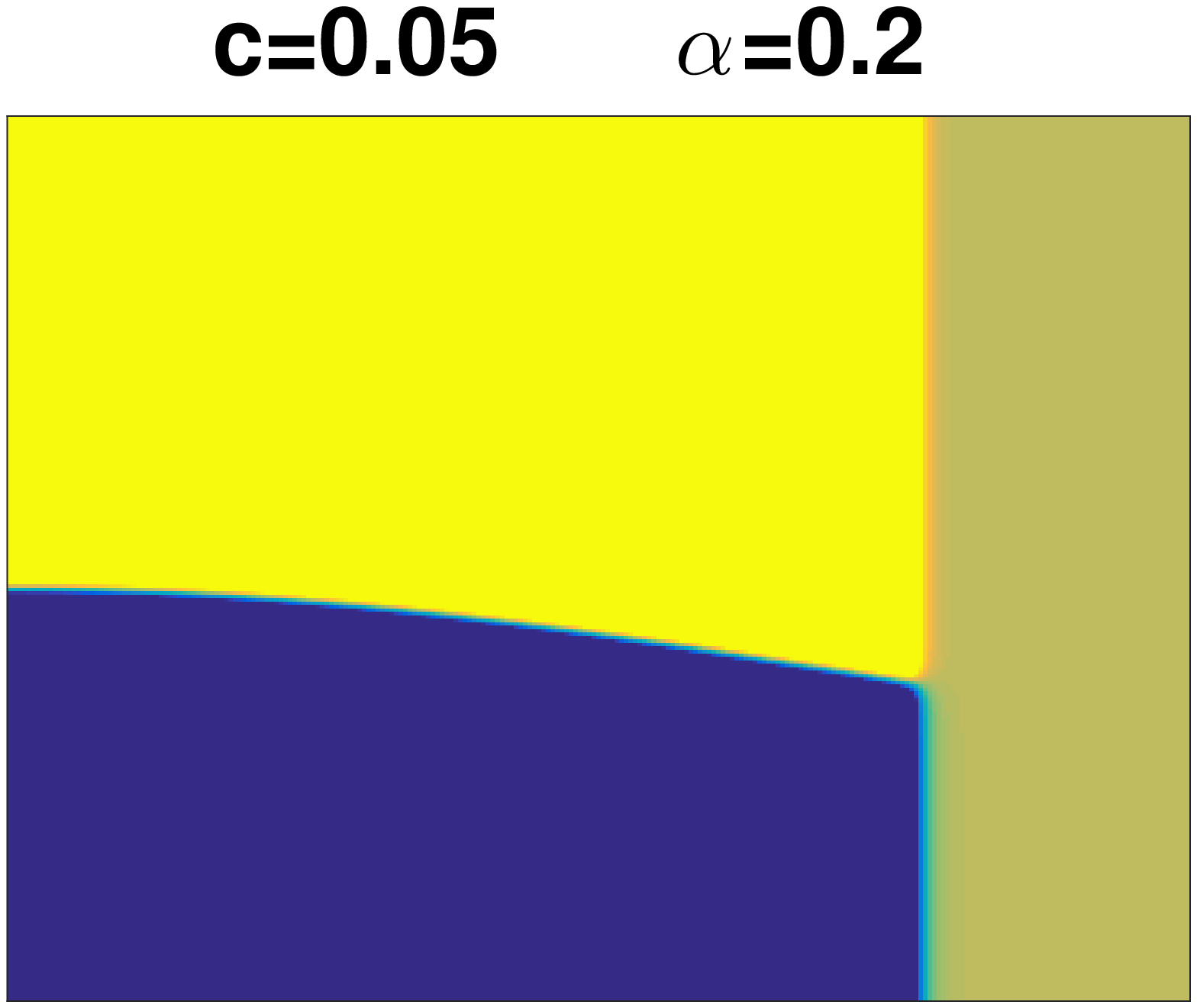}\hfill
\includegraphics[width=0.16\textwidth]{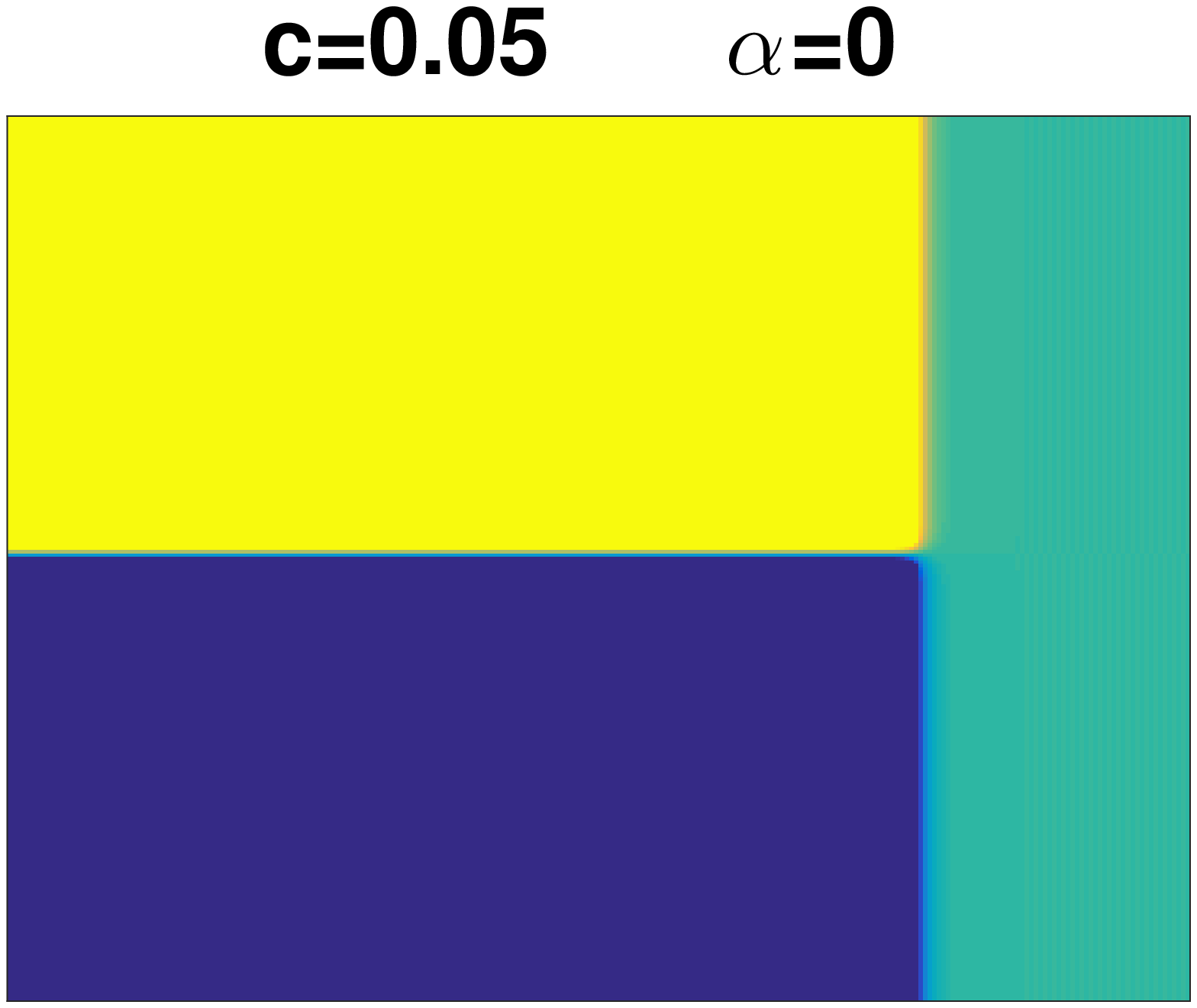}\hfill
\includegraphics[width=0.16\textwidth]{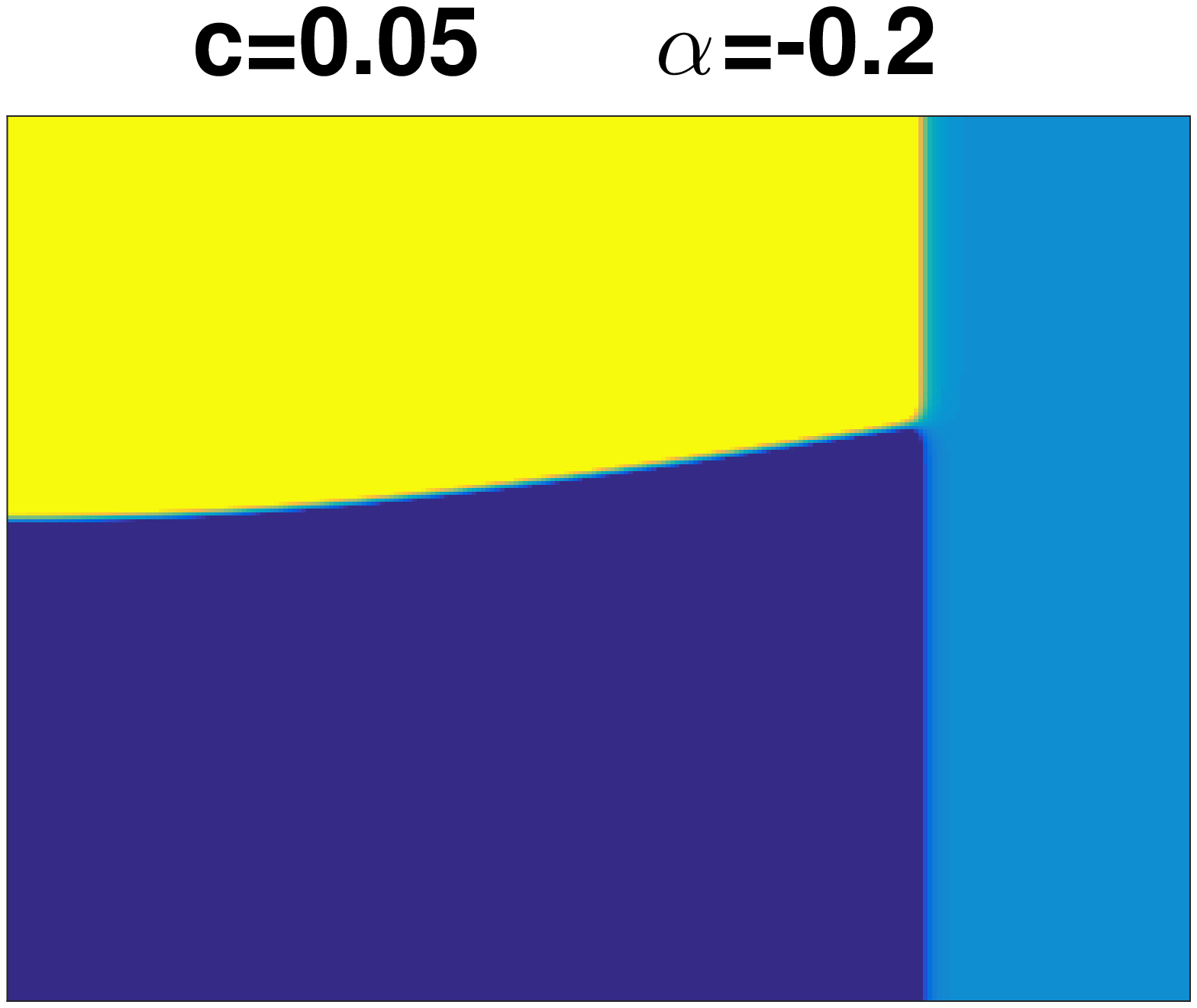}
\caption{Direct simulations with $g(x,u)=1$, $c_x=0.5$, $\alpha=-0.2,0,0.2$ (left) and $c_x=0.5$, $\alpha=-0.2,0,0.2$ (right). Domain is $(0,130)^2$ with Neumann boundary conditions, starting with a step-like initial conditions. }
\label{f:10}
\end{figure}

The result is illustrated in direct simulations in Figure \ref{f:10}, where the selection of the contact angle is clearly visible, dependent on the sign of $\alpha$. One also notices the comparatively weak influence of the left boundary, in particular for large speeds.

We remark here that the restriction on $c_x$ can be relaxed. Our approach is not based on a perturbation from $c_x=0$ but rather exploits an existence result from \cite{Monteiro_Scheel} that was proven for $c_x\gtrsim 0$, only. One expects this existence result to hold for all $c_x<2$ \cite{Monteiro}. We therefore state below a more conceptual result that makes the assumptions required at $\alpha=0$ more explicit. 

There are several  difficulties involved with establishing Theorem \ref{t:1}. First, one notices immediately that the branch of solutions is not continuous in, say, $\mathcal{C}^0(\R^2)$ with the uniform topology, since a changing contact angle will cause values of two different profiles $u$ to differ by roughly 2 in a sector bounded by the two interfaces. A second problem arises when studying the linearization at the solution for $\alpha=0$, which turns out not to be a Fredholm operator. In fact, translation invariance in $y$ implies that the $y$-derivative off this profile belongs to the kernel, but is not localized in space. A classical Weyl sequence argument then implies that the linearization cannot be a Fredholm operator.  Our main technical contributions address these two difficulties as follows. We overcome the lack of Fredholm properties by passing to weighted spaces, relying on a Closed Range Lemma, invertibility in far-field sectors, and a patching argument. We characterize kernels and cokernels exploiting various forms of comparison principles. Related to the first difficulty of a lack of smoothness of the solution is the fact that the resulting Fredholm index is negative. We compensate for the negative index using a farfield-core decomposition, that is, an Ansatz that explicitly inserts a traveling-wave solution in $x<-1$ with a prescribed angle $\varphi$, that eventually compensates for the negative index. 

\paragraph{Outline.} The remainder of this paper is organized as follows. We give a more detailed statement of Theorem \ref{t:1}, making explicit the farfield behavior after formulating conceptual hypotheses in Section \ref{s:2}. Section \ref{s:3} contains the proof of the first main result, establishing the selection of the contact angle based on conceptual assumptions at $\alpha=0$. Section \ref{s:4} establishes the conceptual assumptions in the case of $c_x\ll1$. We conclude with a brief discussion, that also presents some cases where explicit statements for the selected angle are possible.

\section{Traveling waves, farfield asymptotics, and a conceptual perturbation result}\label{s:2}
We first describe far-field asymptotics of solutions by investigating one-dimensional limit equations of \eqref{e:ell} in Section \ref{s:21}. We then state conceptual hypotheses and a more detailed version of Theorem \ref{t:1} in Section \ref{s:22}.

\subsection{Farfield patterns}\label{s:21}
\begin{figure}[htb]
\centering
    \begin{subfigure}[b]{0.52\textwidth}
        \includegraphics[width=\textwidth]{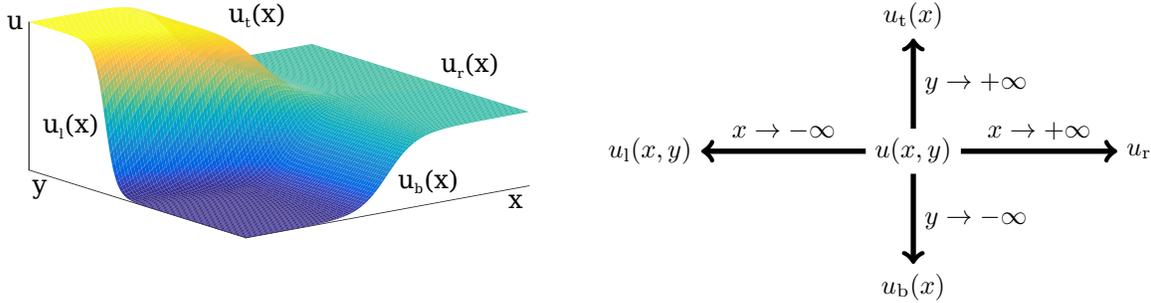}
        \label{fig:Diagram_physical}
    \end{subfigure}
    \begin{subfigure}[b]{0.47\textwidth}
\begin{tikzpicture}[node distance=3cm]
  \node (A) {$u(x,y)$};
  \node [right of=A] (Ri) {$u_\mathrm{r}$};
    \draw[->,line width=2pt] (A) to node[above] {$x\to +\infty$}   (Ri);
  \node [left of=A, xshift=-.5cm] (Le) {$u_\mathrm{l}(x,y)$};
    \draw[->,line width=2pt] (A) to  node[above] {$x\to -\infty$}  (Le);
  \node [above of=A, yshift=-1.2cm] (Up) {$u_\mathrm{t}(x)$};
      \draw[->,line width=2pt] (A) to  node[right] {$y\to +\infty$} (Up);
  \node [below of=A, yshift=1.2cm] (Lo) {$u_\mathrm{b}(x)$};
  \draw[->,line width=2pt] (A) to node[right] {$y\to -\infty$} (Lo);  
\end{tikzpicture}
        \label{fig:mouse}
    \end{subfigure}
\caption{Sketches of the solution with farfield patterns. \label{Fig:H_infty}}
\end{figure}

Recall the definition of the zeros $z_\pm(\alpha)$ and $z_0(\alpha)$, smoothly continuing $\pm1$ and $0$ as zeros for the kinetics in $x<0$ and $x>0$, respectively. Our goal here is to construct solutions to \eqref{e:ell} outside of $\{|x|,|y|\leq 1\}$, that is, in the regions $\{x>1\}$ (right), $\{x<1\}$ (left), $\{y>1\}$ (top), and $\{y<-1\}$ (bottom), for $\alpha\sim 0$; see Figure \ref{Fig:H_infty}.

\paragraph{Farfield pattern to the right.} Here, we simply set $u_\mathrm{r}(\alpha):=z_0(\alpha)$. 

\paragraph{Farfield pattern at the top and bottom.} We seek one-dimensional solutions $u_\mathrm{t}(x;\alpha)$ to \eqref{e:ell} with limits $u_\mathrm{t}(+\infty)=z_0(\alpha)$, $u_\mathrm{t}(-\infty)=z_+(\alpha)$. At $\alpha=0$, such solutions have been constructed in \cite[Prop. 1.4]{Monteiro_Scheel}. In fact, as explained there, the existence follows from a straightforward phase plane analysis for $\alpha=0$, $c_x<2$. The results there also show invertibility of the linearization and monotonicity, $\partial_x u_\mathrm{t}(x;0)<0$. Both properties can be either continued from $\alpha=0$ or established for $\alpha\sim 0$ using the very same methods. In the same way, one establishes the existence of $u_\mathrm{b}(x;\alpha)$, monotonically increasing, with  limits $u_\mathrm{b}(+\infty)=z_0(\alpha)$, $u_\mathrm{b}(-\infty)=z_-(\alpha)$. From a construction using the Implicit Function Theorem in $\mathcal{C}^{1,\beta}$, say, one also finds that both $u_\mathrm{t/b}(x;\alpha)$ are smooth in $\alpha$. Convergence towards the asymptotic states is exponential with uniform rate in $\alpha$, small. 

\paragraph{Farfield pattern to the left.} In $x<-1$, we find the equation 
\begin{equation}\label{e:l}
\Delta u + c_x u_x + c_y u_y+u-u^3+\alpha g_\mathrm{l}(u)=0.
\end{equation}
We look for a family of planar traveling-wave solutions to this equation, $u_\mathrm{l}(x,y;\psi,\alpha)=z_\mathrm{tw}(
\sin(\psi) x + \cos(\psi) y; \psi,\alpha)$, $c_y=c_y(\psi,\alpha)$. Here, the angle $\psi=\varphi-\frac{\pi}{2}$ denotes the deviation of the contact angle from a right angle; see Figure \ref{f:1} and Definition \ref{d:contact}. The profile $z_\mathrm{tw}$ solves the traveling-wave equation,
\begin{equation}\label{e:le}
z_\mathrm{tw}''+(\sin(\psi) c_x+\cos(\psi)c_y)z_\mathrm{tw}'+z_\mathrm{tw}-z_\mathrm{tw}^3+\alpha g_\mathrm{l}(z_\mathrm{tw})=0, \qquad z_\mathrm{tw}(-\infty)=z_-(\alpha), z_\mathrm{tw}(+\infty)=z_+(\alpha).
\end{equation}
Now recall that the ordinary differential equation 
\[
u''+cu'+u-u^3+\alpha g_\mathrm{l}(u)=0,
\]
possess a unique (up to translation) family of solutions with boundary values as in \eqref{e:le}, smoothly depending on $\alpha$ in $\mathcal{C}^2$, with unique $c=c_\mathrm{n}(\alpha)$, smooth. The derivative at $\alpha=0$ can be explicitly calculated as 
\begin{equation}\label{e:c'}
c_\mathrm{n}'(0)=-\int_\R  g_\mathrm{l}(u_*(y))u_*'(y)\rmd y \left(\int_\R \left(u_*'(y)\right)^2\rmd y\right)^{-1},
\end{equation}
with $u_*(y)=\tanh(y/\sqrt{2})$.
We conclude the existence of a solution $z_\mathrm{tw}$ to \eqref{e:le} for all $\psi,\alpha\sim 0$, and 
\begin{equation}\label{e:cy}
c_y(\alpha,\psi)=\frac{1}{\cos(\psi)}c_\mathrm{n}(\alpha)-c_x\tan(\psi).
\end{equation}
One readily finds smooth dependence on $\psi,\alpha$ after an appropriate normalization with respect to translation and exponential convergence in the farfield. 

\paragraph{Patching and partition of unity.}

We define a partition of unity $\chi^j$, $j=\mathrm{t,b,l,r,0}$, that allows us to smoothly glue the solutions constructed here to a solution in the far-field, with small errors. We will work in polar coordinates, $(x,y)=(-r\cos(\psi),r\sin(\psi))$. Note the negative sign in the $x$-variable, corresponding to our orientation of contact angles and interfaces in Figure \ref{f:1} and Definition \ref{d:contact}. Consider the  mollified characteristic functions $\chi_R$ and the partition of unity generated by $\chi_\Psi$, 
\[
\chi_R(r)=\left\{\begin{array}{ll} 1,& r>R\\ 0,&r<R-1\end{array}\right.,\qquad 
\chi_\Psi(\psi)=\displaystyle{\left\{\begin{array}{ll} 1,& \
\frac{\pi}{4}+\frac{1}{10}<\psi<\frac{3\pi}{4}-\frac{1}{10}\\ 0,&\psi>\frac{3\pi}{4}+\frac{1}{10} \mbox{ or } \psi<\frac{\pi}{4}-\frac{1}{10}\end{array}\right.,}
\]
where $\chi_\Psi$ is understood as a smooth function on the circle $\R/2\pi\Z$, and $\sum_{j=1}^4\chi_\Psi(\cdot -j\pi/2)=1$. 

Now construct the farfield partition and the residual in the core through
\begin{align}
\chi_\mathrm{t}(r,\psi)&=\chi_R(r)\cdot \chi_\Psi(\psi),\nonumber\\
\chi_\mathrm{r}(r,\psi)&=\chi_R(r)\cdot \chi_\Psi(\psi-\pi/2),\nonumber\\
\chi_\mathrm{b}(r,\psi)&=\chi_R(r)\cdot \chi_\Psi(\psi-\pi),\nonumber\\
\chi_\mathrm{l}(r,\psi)&=\chi_R(r)\cdot \chi_\Psi(\psi-3\pi/2),\nonumber\\
\chi_0(r,\psi)&=1-\sum_{j=\mathrm{t,r,b,l}}\chi_j(r,\psi).\label{e:pou}
\end{align}
In the following, we will write $\chi_j(x,y)$ for the partition of unity in Cartesian coordinates, slightly abusing notation. We note that $\chi_0$ is compactly supported and that, by 0-homogeneity in the far-field, partial derivatives of the other elements decay algebraically, 
\begin{equation}\label{e:alg}
|\partial^k_x\partial^l_y \chi_j|\leq C (1+r)^{-(k+l)}.
\end{equation}
See also Figure \ref{f:pou} for a sketch. 

In order to state more precise asymptotics for the solution to \eqref{e:ell}, we define the farfield Ansatz as 
\begin{equation}\label{e:ff}
u_\mathrm{ff}(x,y;\psi,\alpha):=\chi_\mathrm{t}(x,y)u_\mathrm{t}(x,y;\alpha)+\chi_\mathrm{r}(x,y)u_\mathrm{r}(x,y;\alpha)+\chi_\mathrm{b}(x,y)u_\mathrm{b}(x,y;\alpha)+\chi_\mathrm{l}(x,y)u_\mathrm{l}(x,y;\psi,\alpha).
\end{equation}
We emphasize that $u_\mathrm{ff}$ is not smooth as an element of $L^\infty(\R^2)$, say, with respect to the parameter $\psi$. It is however smooth in a locally uniform topology, for instance uniform convergence on finite balls $B_R(0)$, as used in Theorem \ref{t:1}.  On the other hand, this farfield Ansatz encodes nodal lines $y\sim -\tan(\psi) x$, hence a contact angle $\varphi=\frac{\pi}{2}+\psi$.

\subsection{Main results and asymptotics}\label{s:22}
We first collect some basic properties of the solution with contact angle $\pi/2$ for $\alpha=0$. We formulate those properties as assumptions. Our first result states that these assumptions hold for $c_x>0$, small.  Our second main result concludes a sharpened version of Theorem \ref{t:1} from these assumptions. 

In order to state hypotheses and main result, we define spaces of exponentially localized functions $L^2_\eta(\R^2)$, and $H^k_\eta(\R^2)$, for any $\eta>0$, as the closure of $\mathcal{C}^\infty_0(\R^2)$ in the norms
\[
\|u(\cdot)\|_{L^2_\eta(\R^2)}:=\|u(\cdot)\rme^{\eta |\cdot |_1}\|_{L^2(\R^2)},\quad \|u(\cdot)\|_{H^k_\eta(\R^2)}:=\sum_{0\leq |\ell|\leq k}\|\partial^\ell u(\cdot)\|_{L^2_\eta(\R^2)},
\]
using multi-index notation $\partial^\ell=\partial_x^{\ell_1}\partial_y^{\ell_2}$, $|\ell|=\ell_1+\ell_1$, and $|(x,y)|_1=|x|+|y|$. 

\begin{Assumption}[Existence and properties of solutions for $\alpha=0$, $c_x>0$]\label{h:1}
Fix $c_x>0$ and set $\alpha=0$. 
\begin{itemize}
 \item [(A1)] \emph{Existence:} We assume that there exists a solution $\Theta(x,y)$ to  \eqref{e:ell}  with $c_y=0$, $\varphi=\pi/2$; moreover, assume that for any fixed $x$ the mapping $y \mapsto \Theta(x,y)$ is nondecreasing.
 \item [(A2)] \emph{Asymptotics:} We assume that $\Theta(x,y)$ converges exponentially to asymptotics profiles $u_\mathrm{r}(y;0)$, $u_\mathrm{l}(y;0,0)$ and $u_\mathrm{t,b}(x;0)$, for $x\to\pm\infty$ and $y\to\pm\infty$, respectively. More precisely, 
$(\Theta-u_j)\chi_j \in H^2_\eta(\R^2)$ for $j=\mathrm{t,b,l,r}$ and some $\eta>0$.
  \item [(A3)]  \emph{Linearization:} The operator 
 \[
\mathscr{L}: H_\eta^2(\R^2)  \to L^2_\eta(\R^2), \quad u\mapsto \Delta u + c_x u + \mu(x)u-3\Theta^2u.
\]
is  Fredholm with index $-1$ for all $\eta>0$, sufficiently small, with trivial kernel, and with cokernel spanned by $e^{c_x x}\Theta_y(\cdot,\cdot)\in L^2(\R^2)$.
\end{itemize}
\end{Assumption}

\begin{Theorem}[Existence of balanced fronts]\label{Theorem:small_c_x} Assumptions  (A1), (A2), and (A3) hold for all  $c_x > 0$ sufficiently small. 
\end{Theorem}

It turns out that Assumptions (A1) and (A2) are rather direct consequences of the results in \cite{Monteiro_Scheel}, for small $c_x>0$. Establishing (A3) will take up the major part of Section \ref{s:4}.

\begin{Theorem}[Existence of oblique fronts --- refined asymptotics]\label{Theorem:all_c_x}
Assume (A1), (A2), and (A3), for some fixed $c_x>0$. Then there exists $\alpha_0>0$ such that there is a solution  $u_*(x,y;\alpha)$ to \eqref{e:ell}, with contact angle  $\varphi(\alpha)$ as in Definition \ref{d:contact} and $c_y=c_y(\alpha,\varphi(\alpha)-\frac{\pi}{2})$ as in \eqref{e:cy}, for all $|\alpha|\leq \alpha_0$. 
More precisely, we have that, writing $\psi(\alpha)=\varphi(\alpha)-\frac{\pi}{2}$,
\begin{enumerate}
\item $\varphi(\alpha)$ is smooth with $\varphi(0)=\pi/2$;
\item $u_*(x,y;\alpha)=u_\mathrm{ff}(x,y;\alpha,
\psi(\alpha))+w(x,y;\alpha)$, where $u_\mathrm{ff}$ was defined in \eqref{e:ff} and $\alpha\mapsto w(\cdot;\alpha)\in H^2_\eta(\R^2)$ is smooth, well-defined for $\eta>0$ sufficiently small.
\end{enumerate}
\end{Theorem}
It is straightforward to verify that Theorem \ref{Theorem:all_c_x} implies Theorem \ref{t:1}.

\section{Angle selection as a perturbation result}\label{s:3}
We proof Theorem \ref{Theorem:all_c_x} in this section. We first introduce a shear coordinate transformation such that $u_\mathrm{l}$ is independent of $x$ in the farfield. We then set up an implicit function theorem and show that the nonlinear mapping is well-defined. The key step to applying the implicit function theorem consists of establishing invertibility of the linearization. Since the linearization $\mathscr{L}$ is Fredholm of index $-1$, we use a bordering lemma, adding $\varphi$ as an additional variable, to obtain invertibility and conclude Theorem \ref{Theorem:all_c_x}.

\paragraph{Shear coordinates.}
 We define a shear transformation:
\begin{equation}\label{partition_of_unity_for_shear_transformation}
(\tx,\ty) = \mathbb{S}^{(\psi)}(x,y) = \left(x,  y + x \chi^-(x)\tan(\psi)\right).
\end{equation}
where $\chi^-(x)=1, x<-2$ and  $\chi^-(x)=0, x>-1$, is a smooth function. Note that  $\mathbb{S}^{(0)}(\cdot, \cdot)$ is the identity map and, inverting explicitly,  $\mathbb{S}^{(\psi)}(\cdot,\cdot)$ is a $\mathcal{C}^{\infty}$-diffeomorphism of $\mathbb{R}^2$, for each fixed $\psi \in \left(-\frac{\pi}{2}, \frac{\pi}{2} \right)$.

Straightforward calculations show that, in the new coordinates, \eqref{e:ell} becomes
\begin{align}\label{e:ellt}
\left[\Delta_{\tx,\ty} + c_x \partial_{\tx} + c_y \partial_{\ty} + \mu(x)\right]\tilde{u} - \tilde{u}^3  + \alpha  g(\tx,\tilde{u}) + \tan(\psi)\left[ 2 \mathcal{S}_x \partial_{\tx \, \ty} + (\mathcal{S}_{x})^2\tan(\psi)\partial_{\ty\ty} +  (c_x \mathcal{S}_x + \mathcal{S}_{xx})\partial_{\ty}\right]\tilde{u} =0,
\end{align}
where $\mathcal{S}(x) = x \chi^-(x)$ has linear growth, although all its derivatives are bounded. At  $\psi = 0$, we recover \eqref{e:ell}.

\begin{figure}
\hspace*{0.8in}\includegraphics[width=0.14\textwidth]{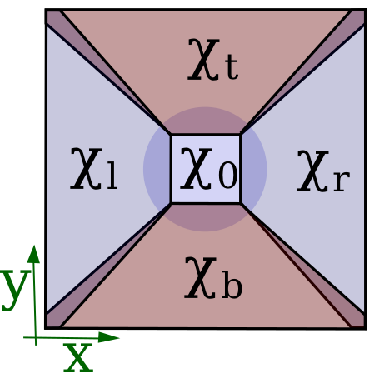}\qquad
\includegraphics[width=0.14\textwidth]{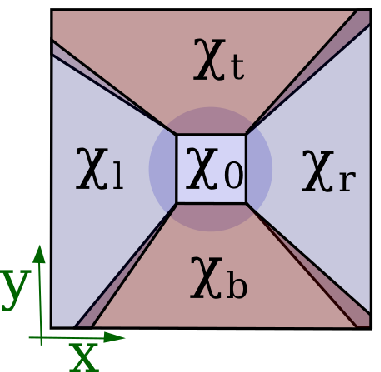}\hfill
\raisebox{0.15in}{\includegraphics[width=0.24\textwidth]{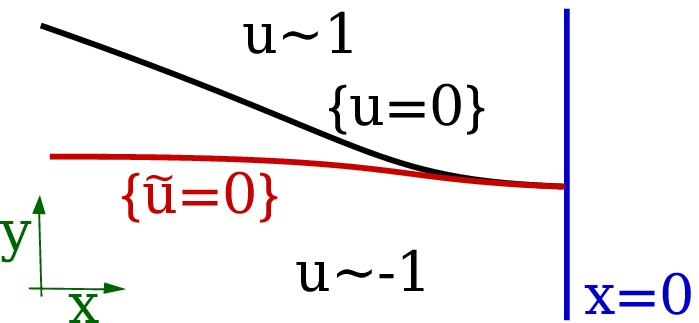}}\hspace*{0.8in}
\caption{Partition of unity (left) and sheared partition of unity (middle), as well as interfacial lines before and after shearing (right).}
\label{f:pou}
\end{figure}


\paragraph{Farfield-core decomposition and smoothness of the nonlinear mapping.}

We will substitute an Ansatz $u={u}_\mathrm{ff}+w$ into \eqref{e:ell}, consider the resulting equation in the sheared coordinates,  and strive to solve for $(w,\psi)$, as functions of $\alpha$, exploiting the choice of $c_y$ from \eqref{e:cy} for compatibility in the region $x<-1$. 

In the transformed coordinates, we obtain a new partition of unity $\tilde{\chi}_\mathrm{t,r,b,l,0}$, depending on the angle $\psi$ in a mild fashion, that is, derivatives in $\tilde{x},\tilde{y}$ still exhibit decay as stated in \eqref{e:alg} and derivatives with respect to $\psi$ are bounded. The blending regions along angles $\psi\sim \pi/4$ and $\psi\sim -\pi/4$  are slightly realigned; see Figure \ref{f:pou}.

More importantly, the Ansatz functions $u_\mathrm{t/b/r}$ are unaffected by the shear coordinate change. The Ansatz function $u_\mathrm{l}$ simplifies to $u_\mathrm{l}(\cos(\psi)\tilde{y};\psi,\alpha)$, upon using the same angle $\psi$ in the Ansatz for $u_\mathrm{l}$ and for the shear transformation. As a consequence, exploiting uniform exponential decay of derivatives of $u_\mathrm{l}$, the transformed function  $\tilde{u}_\mathrm{l}$ is smoothly dependent on $\psi$ in $\mathcal{C}^2$. 

We denote by $\tilde{u}_\mathrm{ff}$ the farfield correction in the shear coordinates $(\tilde{x},\tilde{y})$. and write $\tilde{w}$ for the core correction. Substituting $u=\tilde{u}_\mathrm{ff}+\tilde{w}$ into \eqref{e:ellt} and setting $c_y(\alpha,\psi)$ as in \eqref{e:cy} gives a nonlinear equation with variables $(w,\psi,\alpha)$, which reads, after dropping tildes for readability,
\begin{align}\label{e:ellta}
0=&\left[\Delta + c_x \partial_{x} + c_y \partial_{y} + \mu(x)\right]({w}+{u}_\mathrm{ff}) - ({w}+{u}_\mathrm{ff})^3  \nonumber\\
&+ \alpha  g(x,{w}+{u}_\mathrm{ff}) + \tan(\psi)\left[ 2 \mathcal{S}_x \partial_{xy}^2 + (\mathcal{S}_{x})^2\tan(\psi)\partial_{y}^2 +  (c_x \mathcal{S}_x + \mathcal{S}_{xx})\partial_{y}\right]({w}+{u}_\mathrm{ff})\nonumber \\
=:&F(w,\psi,\alpha),
\end{align}
where $u_\mathrm{ff}$ depends implicitly on $\alpha,\psi$. By assumption (A1), 
$w=w_0:=\Theta-u_\mathrm{ff}$ is a solution at $\psi=\alpha=0$. 

\begin{Lemma}[Smoothness]\label{l:sm}
The function $F$ defined in \eqref{e:ellta} is smooth as a mapping from a neighborhood of $(w_0,0,0)$ in $H^2_\eta(\R^2)\times\R^2$ into $L^2_\eta(\R^2)$ for $\eta>0$, sufficiently small.
\end{Lemma}
\begin{Proof}
We first show that $F$ is well defined. Clearly, $F$ belongs to $L^2_\mathrm{loc}(\R^2)$ when $w\in H^2(\R^2)$. We show that $F(w_0,\psi,\alpha)\in L^2_\eta(\R^2)$.  This is a consequence of the construction of farfield profiles in our Ansatz as we will see next. In the regions where the partition of unity does not overlap, we have that  $u_\mathrm{ff}$ is an exact solution and $F(w_0,\psi,\alpha)$ therefore vanishes. Due to exponential convergence of the Ansatz functions $u_\mathrm{t,r,b,l}$ to their asymptotic states, errors in the overlapping regions are exponentially small.  Next, note that $w_0\in H^2_\eta(\R^2)$ due to (A2). As a consequence,  $F(w_0,\psi,\alpha)\in L^2_\eta(\R^2)$. 

The  linear terms in $w$ belong to $L^2_\eta(\R^2)$ provided that $(x,y)$-dependent coefficients are bounded. This in turn is a consequence of the fact that derivatives of $\mathcal{S}$ are bounded. Nonlinear terms involving $w$ are automatically in $L^2_\eta(\R^2)$ since $H^2_\eta(\R^2)$ is an algebra. 

Continuity in $w,\psi,\alpha$ follows from the smooth dependence in the equation and uniform exponential bounds in the farfield. One similarly obtains that derivatives with respect to  $w,\psi,\alpha$ are well-defined, bounded, and continuous in $w,\psi,\alpha$.
\end{Proof}

We remark that it is here that we exploit the shear transformation. In the original coordinates, a $\psi$-derivative would generate a term $xw$ in $x<-2$, which is unbounded in $L^2_\eta(\R^2)$. 

\begin{Lemma}[Invertibility of the linearization]\label{l:inv}
The linearization $D_{w,\psi}F|_{(w_0,0,0)}:H^2_\eta(\R^2)\times\R\to L^2_\eta(\R^2)$ is bounded invertible.
\end{Lemma}
\begin{Proof}
By (A3), $D_wF(w_0,0,0)=\mathscr{L}$ is Fredholm of index -1, with trivial kernel, such that it is sufficient to show that $D_\psi F(w_0,0,0)$ does not belong to the range of $D_wF(w_0,0,0)$. Using the explicit expression for the cokernel from (A3), we find
\[
D_{\psi}F|_{(w_0,0,0)}=\partial_wF(w_0,0,0)\partial_\psi u_\mathrm{ff}+\left(2 \mathcal{S}_x \partial_{x y} +  (c_x \mathcal{S}_x + \mathcal{S}_{xx})\partial_{y}\right)\Theta-c_x\Theta_y,
\]
where the last term stems from differentiating $c_y\Theta_y$, with explicit derivative $\partial_\psi c_y=-c_x$ at $\psi=\alpha=0$ from \eqref{e:cy}. 
Again by (A3), the cokernel of $\mathscr{L}$ is spanned by $\Theta_y\rme^{c_x x}$. We therefore need to evaluate the $L^2$-inner product between $D_\psi F(w_0,0,0)$ and $\Theta_y\rme^{c_x x}$,
\begin{align*}
M_\psi=&\int_{\R^2} \Theta_y \rme^{c_x x} \left\{\mathscr{L}(\partial_\psi u_\mathrm{ff})+2 \mathcal{S}_x \Theta_{xy}  +  (c_x \mathcal{S}_x + \mathcal{S}_{xx})\Theta_y-c_x\Theta_y\right\}\rmd x\rmd y\\
=&\int_{\R^2} \Theta_y \rme^{c_x x} \left\{2 \mathcal{S}_x \Theta_{xy} +  (c_x \mathcal{S}_x + \mathcal{S}_{xx})\partial_{y}\Theta-c_x\Theta_y\right\}\rmd x\rmd y.
\end{align*}
Here we used that $\partial_\psi u_\mathrm{ff}$ is exponentially localized, belongs in particular to $H^2_\eta(\R^2)$ for $\eta>0$ sufficiently small, and the scalar product of $\mathscr{L}(\partial_\psi u_\mathrm{ff})$ with the cokernel vanishes as a consequence. Writing $2\Theta_{xy} \Theta_y=\partial_x\Theta_y^2$ and integrating by parts, exploiting that boundary terms vanish due to the exponential factor, we find 
\begin{equation}\label{e:mel}
M_\psi=-c_x\int_{\R^2}\left(\Theta_y\right)^2\rme^{c_x x}\rmd x \rmd y<0.
\end{equation}
This proves the lemma. 
\end{Proof}
We are now ready to proof our first main theorem. 
\begin{Proof}[of Theorem \ref{Theorem:all_c_x}]
Using Lemma \ref{l:sm} and \ref{l:inv}, we can use the Implicit Function Theorem to solve $F(w,\psi,\alpha)=0$ and find a branch of solutions $(w,\psi)(\alpha)$. The decomposition stated in the theorem is an immediate consequence of the farfield-core decomposition used in the proof. 
\end{Proof}
On a coarse level, the most interesting information here is of course the contact angle $\varphi=\frac{\pi}{2}+\psi$. Its derivative with respect to the perturbation parameter $\alpha$ can be readily obtained by projecting leading-order terms on the cokernel, 
\begin{equation}\label{e:dp}
\frac{\rmd\varphi}{\rmd\alpha}=-\frac{M_\alpha}{M_\psi},\qquad M_\alpha=\int_{\R^2} \partial_\alpha F \Theta_y \rme^{c_x x} \rmd x \rmd y,
\end{equation}
evaluated at $(w_0,0,0)$. For applications in Section \ref{s:5.1}, we compute 
\[
\partial_\alpha F=\mathscr{L}(\partial_\alpha) u_\mathrm{ff}+\partial_\alpha c_y \Theta_y + g(x,\Theta).
\]
Using the geometric relation \eqref{e:cy}, and recalling the definition of the normal speed $c_\mathrm{n}(\alpha)$, we find $\partial_\alpha c_y=c_\mathrm{n}'$ at $\alpha=\psi=0$. Further exploiting that $\partial_\alpha u_\mathrm{ff}$ belongs to $H^2_\eta(\R^2)$, we find after a short calculation
\begin{equation}\label{e:ma}
M_\alpha=-c_\mathrm{n}'(0)M_\psi-\int_\R \rme^{c_x x} \left[G(x,u_\mathrm{t})-G(x,u_\mathrm{b})\right] \rmd x,
\end{equation}
where $\partial_u G(x,u)=-g(x,u)=\partial_\alpha W(x,u)$ and $c_\mathrm{n}'(0)$ is given in \eqref{e:c'}.

\section{Establishing (A1)--(A3) and Fredholm properties of the linearization}\label{s:4}
In this section, we proof Theorem \ref{Theorem:small_c_x}, that is, we establish assumptions (A1)--(A3) for small speeds, $c_x>0$. In particular, throughout this section, we will work with $\alpha=0$. We first recall the main relevant result from \cite{Monteiro_Scheel}, Section \ref{s:3.0}, which establishes (A1) and gives some additional qualitative properties.  We show that $\mathscr{L}$ is Fredholm in $L^2_\eta(\R^2)$ in Section \ref{s:3.1}, and we establish exponential asymptotics (A2), in Section \ref{s:3.2}. Finally, we compute the Fredholm index and show that $\mathscr{L}$ has trivial kernel in Section \ref{s:3.3}, establishing (A3). 

\subsection{Existence and qualitative properties at small speeds}\label{s:3.0}
We recall the relevant parts of \cite[Prop. 1.4]{Monteiro_Scheel} and prove some slight refinements.

\begin{Proposition}\label{p:0}
For all $c_x\geq 0$, sufficiently small, there exists a solution $\Theta$ to \eqref{e:ell} at $\alpha=0$ with contact angle $\pi/2$. In addition, we have that $\Theta(x,y)=-\Theta(x,-y)$, $\Theta_y(x,y)>0$ for all $(x,y)$, and 
\begin{align*}
\lim_{x\to -\infty}\Theta(x,y)&=u_\mathrm{l}(y;0,0)=\tan(y/\sqrt{2}), \mbox{ uniformly in }y;\\
\lim_{x\to \infty}\Theta(x,y)&=u_\mathrm{r}(0)=0, \mbox{ uniformly in }y;\\
\lim_{y\to -\infty}\Theta(x,y)&=u_\mathrm{b}(x;0), \mbox{ uniformly in }x;\\
\lim_{y\to -\infty}\Theta(x,y)&=u_\mathrm{t}(x;0,0), \mbox{ uniformly in }x.
\end{align*}
\end{Proposition}
\begin{Proof}
Existence, reflection, and convergence properties have been established in \cite{Monteiro_Scheel}, as well as weak monotonicity $\Theta_y(x,y)\geq0$ at $c_x=0$. To show strong monotonicity for $c_x>0$, first notice that weak monotonicity implies strong monotonicity as follows.  We use a  Harnack type inequality \cite[Thm. 9.22]{gilbarg2015elliptic}. Suppose $u$ is nonnegative, $Lu \leq f $, for $u, f \in W^{2,n}(\Omega)$ with  $B_{2r} \subset \Omega$, then
\begin{equation*}
\left(\frac{1}{|B_r|}\int_{B_r}u^p\right)^{1/p} \leq C \left(\inf_{B_r}u +r \| f\|_{B_{2r}}\right) .
\end{equation*} 
Applied  locally with $f=0$, $u=\Theta_y$ and $L = \mathscr{L}$, this implies that  $\{(x,y)|\Theta_y(x,y) =0\}$ is open. Since this set is clearly also closed, it is empty as a subset of the (connected) plane $\R^2$, hence $\Theta_y(x,y)>0$ as claimed.

It remains to show that $\Theta_y\geq 0$ for $c_x>0$, small enough. First, for any $R>0$, there is $c_0>0$ such that $\Theta_y>0$ in $\Omega=\{(x,y|\,|y|\leq R,\ x<R\}$, by positivity at $c_x=0$ and continuity as established in \cite{Monteiro_Scheel}. Also, note that, for all $c_x\geq 0$, small, we have  $\liminf_{|(x,y)|\to\infty} \Theta_y \geq 0$. We next construct a function $\bar{w}(x)>0$ that will serve as a supersolution in the complement $\Omega^\mathrm{c}$. We want to solve 
\[
w_{xx}+c_xw_x+(\mu-3u_\mathrm{t}^2)w=-1.
\]
The left-hand side defines a Sturm-Liouville operator that we claim possesses strictly negative spectrum. One easily finds that the essential spectrum has negative real part considering the limnits $x=\pm\infty$, and that point spectrum is real. We showed in \cite{Monteiro_Scheel} that the linearization is negative definite at $c_x=0$. The construction of profiles immediately gives continuity of profiles in $L^\infty(\R)$ as a function of $c_x$, such that point spectra are continuous in $c_x$. On the other hand, the construction also shows that $\lambda=0$ is not an eigenvalue for any $c_x$. This can be seen by inspecting the phase portrait for finite $c_x$, where stable and unstable manifolds intersect necessarily in a transverse fashion, thus implying that the linearization does not possess a bounded solution.  We conclude that the spectrum has strictly negative real part for all $c_x\geq 0$ since eigenvalues are all negative for $c_x=0$ and cannot cross $\lambda=\rmi\R$ for $c_x>0$ since they are real and nonzero, as explained above. 
As a consequence, the Green's function is negative and the solution $w$ is positive, approaching  nonzero limits $w_\pm>0$ at $x=\pm\infty$ exponentially. We smoothly change $w$ to $\bar{w}$ such that  $\bar{w}\equiv w_+$ in $x>R$, while preserving the fact that $w$ is a supersolution
\[
\bar{w}_{xx}+c_x\bar{w}_x+(\mu-3u_\mathrm{t}^2)\bar{w}\leq -1/2.
\]
Choosing $R$ large enough and exploiting boundedness of $\bar{w}$ as well as uniform convergence $\Theta(x,y)\to u_\mathrm{t}(x)$, we also have 
\[
\Delta\bar{w}+c_x\bar{w}_x+(\mu-3\Theta^2)\bar{w}\leq -1/4,\qquad \text{ in } |y|\geq R.
\]
One also readily finds that 
\begin{equation}\label{e:bw}
\Delta\bar{w}+c_x\bar{w}_x+(\mu-3\Theta^2)\bar{w}\leq -1/4,\qquad \text{ in } x\geq R,
\end{equation}
thus establishing that $\bar{w}$ is a supersolution in $\Omega^\mathrm{c}$. Also, we have $\liminf_{|(x,y)|\to\infty} \bar{w} \geq \delta>0$, for some $\delta$ independent of $c_x$ and $R$. Now, consider $v:=\Theta_y+\kappa\bar{w}$. Define 
\[
\kappa_*=\inf\{\kappa\geq 0|v>0\}. 
\]
Clearly, the infimum is taken over a nonempty set and is therefore well-defined. If $\kappa_*=0$, then $\Theta_y\geq 0$. We therefore assume $\kappa_*>0$. Using the fact that $v>0$ in $\Omega$ and that $v>0$ at infinity, we find a zero minimum at $(x_*,y_*)\in \Omega^\mathrm{c}$, which implies that, at $(x,y)=(x_*,y_*)$,
\[
\Delta v+c_xv_x+(\mu-3\Theta^2)v\geq 0.
\]
On the other hand, since $\Theta_y$ is in the kernel,
\[
\Delta v+c_xv_x+(\mu-3\Theta^2)v= 
\kappa\left(\Delta \bar{w}+c_x\bar{w}_x+(\mu-3\Theta^2)\bar{w}\right)\leq -\kappa/4,
\]
in $\Omega^\mathrm{c}$, establishing a contradiction. Hence, $\kappa_*=0$ and $\Theta_y\geq 0$ as claimed. 
\end{Proof}

\subsection{The linearization is Fredholm}\label{s:3.1}
We rely on the following abstract result to prove that the linearization is Fredholm in appropriately weighted spaces.  
\begin{Lemma}[Closed Range Lemma {\cite[Prop. 6.7]{Taylor}}]\label{Abstract_Closed_Range_Lemma}
Given a sequence of Banach spaces $X \subset Y \subset Z$ so that $X \hookrightarrow Y$ is continuous and dense and  $X \hookrightarrow Z$ is continuous and dense, let $K: X \to Z$ be a compact linear operator and $T: X \to Y$ be a continuous linear operator. If 
\begin{equation*}  
\|u\|_X \leq C\|Ku \|_Z + \|Tu\|_Y,
\end{equation*}
then $T$ is a semi-Fredholm operator, i.e., it has closed range and finite dimensional kernel. 
\end{Lemma}
In our case, the compactness portion stems from contributions in a bounded region of the plane. We show next how to apply this lemma in a somewhat more general case of an elliptic operator with coefficients that have limits as $x,y\to\pm\infty$. Consider therefore 
\begin{equation}\label{e:a}
\mathcal{A} := \Delta + b(x,y)\cdot \nabla  + a(x,y),
\end{equation}
with $a,b\in L^\infty(\R^2)$. We assume that $a,b$ possess uniform limits,
\begin{equation}\label{e:al}\begin{array}{llll}
\lim_{x\to -\infty} a(x,y)=a_\mathrm{l}(y),& \lim_{x\to \infty} a(x,y)=a_\mathrm{r}(y), & 
\lim_{y\to -\infty} a(x,y)=a_\mathrm{b}(x),& \lim_{y\to \infty} a(x,y)=a_\mathrm{t}(x),\\
\lim_{x\to -\infty} b(x,y)=b_\mathrm{l}(y),& \lim_{x\to \infty} b(x,y)=b_\mathrm{r}(y),& 
\lim_{y\to -\infty} b(x,y)=b_\mathrm{b}(x),& \lim_{y\to \infty} b(x,y)=b_\mathrm{t}(x),
\end{array}
\end{equation}
and define the limiting operators 
\begin{equation}\label{e:alim}
\mathcal{A}_ju=\Delta u + b_j\cdot \nabla u  + a_j u,\qquad j=\mathrm{l,r,b,t}.
\end{equation}
%
%
%
%
\begin{Proposition}[Asymptotic invertibility implies Fredholm]\label{Proposition:closed_range_elliptic_regularity}
Assume that $\mathcal{A}$, defined in \eqref{e:a} with domain of definition $H^2(\R^2)$ in $L^2(\R^2)$, possesses limits as in \eqref{e:al}, and that the limiting operators $\mathcal{A}_j$ from \eqref{e:alim} are bounded invertible in $L^2(\R^2)$.  Then $\mathcal{A}$ is semi-Fredholm. In particular, we have, for any $R>0$ sufficiently large, that there is a constant $C(R)$ such that 
\begin{equation*}
\|u \|_{H^2(\mathbb{R}^2)}\leq C(R)\left( \| u\|_{L^2(B_R)} + \|\mathcal{A}u \|_{L^2(\mathbb{R}^2)}\right),
\end{equation*}
where $B_R:=\{|x|\leq R\}.$
\end{Proposition}

\begin{Proof}
Here, $C$ denotes a constant that may change throughout but does not depend on quantities appearing in the equation unless otherwise noted. Elliptic regularity readily gives 
\begin{equation}\label{elliptic_regularity_inequality}
\|u \|_{H^2(\mathbb{R}^2)} \leq C( \| u\|_{L^2(\mathbb{R}^2)} + \|\mathcal{A}u \|_{L^2(\mathbb{R}^2)}),
\end{equation}
for some $C>0$. In the following, we briefly write  $\|\cdot\|_{H^2}$ for $\|\cdot\|_{H^2(\mathbb{R}^2)}$, and similarly for $L^2$ and $H^1$ norms. We begin by splitting the first term on the right-hand side 
\[
\| u\|_{L^2} \leq \| u\|_{L^2(B_R)} + \sum_{j= \mathrm{l,r,b,t}}\|\chi_j u\|_{L^2},
\]
where the $\chi_j$ are elements of the partition of unity \eqref{e:pou}, supported in $|(x,y)|\geq R-1$. Next, decomposing $\mathcal{A}u = f$ in the equivalent form $\mathcal{A} \chi_j u = \chi_j f + [\mathcal{A}, \chi_j]u$, where brackets denote the commutator,  we have that 
\begin{align*}
\mathcal{A}_j(\chi_ju) &= \chi_jf + [\mathcal{A}, \chi_j]u  +  [\mathcal{A}_j - \mathcal{A}]\chi_ju.
\end{align*}
At this point we use that the far field operators $\mathcal{A}_j$ are bounded invertible, obtaining
\begin{align*}
\|\chi_ju\|_{H^2} &\leq \|\left(\mathcal{A}_j\right)^{-1}\chi_jf\|_{H^2} + \|\left(\mathcal{A}_j\right)^{-1}[\mathcal{A}, \chi_j]u\|_{H^2}  +  \|\left(\mathcal{A}_j\right)^{-1}[\mathcal{A}_j - \mathcal{A}](\chi_ju)\|_{H^2} \nonumber \\
& \leq C\left(\|\chi_j f\|_{L^2} +\|[\mathcal{A}, \chi_j]u\||_{L^2} + \|[\mathcal{A}_j - \mathcal{A}]\chi_ju\|_{L^2} \right).
\end{align*}
By convergence, we can choose $R$ sufficiently large such that $\|(\mathcal{A}-\mathcal{A}_j)\chi_j\|_{H^1\to L^2}\leq \varepsilon$, arbitrarily small.
Hence,  
\[
 \|\chi_j u\|_{H^2} \leq C \left( \|\chi_j f\|_{L^2} +\|[\mathcal{A}, \chi_j]u\||_{L^2}\right).
 \]
Inserting this estimate into  \eqref{elliptic_regularity_inequality}, we obtain
\begin{align}\label{almost_final_equation_elliptic_regularity_inequ}
\|u \|_{H^2} &\leq C\left( \| u\|_{L^2(B_R)} + \sum_{d=1}^4\|\chi_j u\|_{L^2} + \|\mathcal{A}u \|_{L^2} \right)\nonumber\\
&\leq  C\left(\| u\|_{L^2(B_R)} + \|\mathcal{A}u \|_{L^2} + \sum_j \|\chi_j f\|_{L^2} +\|[\mathcal{A}, \chi_j]u\||_{L^2}\right) \nonumber\\
& \leq C\left( \| u\|_{L^2(B_R)} + \|\mathcal{A}u \|_{L^2}  +\|[\mathcal{A}, \chi_j]u\||_{L^2}\right).
\end{align}
Using smallness of derivatives of $\chi_j$ when $R$ is large \eqref{e:alg}, we find 
\begin{align*}
\|[\mathcal{A}, \chi_j]u\||_{L^2} \leq \sum_{1\leq |\alpha |\leq 2, j =1,\ldots 4,} \|\partial^{\alpha }\chi_j\|_{L^{\infty}}\cdot\|u\|_{H^1(\mathbb{R}^2)} \leq \varepsilon \|u\|_{H^1},  
\end{align*}
with $\varepsilon$ arbitrarily small when $R$ is large. Absorbing this term on the left hand side of \eqref{almost_final_equation_elliptic_regularity_inequ}, we obtain  $\displaystyle{\|u \|_{H^2(\mathbb{R}^2)}  \leq C( \| u\|_{L^2(B_R)} + \|\mathcal{A}u \|_{L^2(\mathbb{R}^2)}}),$ as claimed. Together with Lemma \ref{Abstract_Closed_Range_Lemma}, this establishes that $\mathcal{A}$ is semi-Fredholm. 
\end{Proof}
\begin{Remark}
From the proof, we see that it is sufficient to require the slightly weaker convergence $\|(a-a_j)\chi_j\|_{L^\infty(\R^2)}\to 0$, $\|(b-b_j)\chi_j\|_{L^\infty(\R^2)}\to 0$  as the perimeter $R$ of the partition of unity tends to infinity.
\end{Remark}
\begin{Corollary}\label{c:1}
Under the assumption of Proposition \ref{Proposition:closed_range_elliptic_regularity}, suppose that $b,\nabla b\in L^\infty(\R^2)$ satisfy a convergence estimate of the form \eqref{e:al}. Then $\mathcal{A}$ is Fredholm. 
\end{Corollary}
\begin{Proof}
We apply Proposition \ref{Proposition:closed_range_elliptic_regularity} to the $L^2$-adjoint $\mathcal{A}^*u=\Delta u + \nabla\cdot(b u)+au$, which is of the same form as $\mathcal{A}$, with limiting operators simply being the adjoints of the limiting operators $\mathcal{A}_j$. As a conequence, $\mathcal{A}^*$ is semi-Fredholm, and hence the cokernel of $\mathcal{A}$ is finite-dimensional, establishing that $\mathcal{A}$ is Fredholm as claimed. 
\end{Proof}

%

Unfortunately, Corollary \ref{c:1} cannot be directly applied to $\mathscr{L}$, since $\mathcal{A}_\mathrm{l}$ is not invertibility due to a kernel, spanned by $\Theta_y(x=-\infty,y)$.  We therefore resort to exponentially weighted spaces that allow for control of localization of functions along the vertical quenching line and along the interface $y\sim 0$, respectively. It is convenient to slightly generalize the class of exponential weights considered. Consider therefore the smooth rate functions 
\[
\sigma_{\eta_\mathrm{l},\eta_\mathrm{r}}(x)=-\eta_\mathrm{l}x\chi_-(x)+\eta_\mathrm{r}x\chi_+(x),\qquad
\sigma_{\eta_\mathrm{b},\eta_\mathrm{t}}(y)=-\eta_\mathrm{b}y\chi_-(y)+\eta_\mathrm{t}y\chi_+(y),
\]
with $\chi_\pm$ a smooth partition of unity for the real line mollifying the indicator functions of $\R_\pm$. Define the associated exponentially weighted spaces $L^2_{{\underline{\eta}}}(\R^2)$, ${\underline{\eta}}=(\eta_\mathrm{l},\eta_\mathrm{t},\eta_\mathrm{b},\eta_\mathrm{t})$, with norm
\[
\|u\|^2_{L^2_{\underline{\eta}}}=\int_{\R^2} \left|u(x,y)\rho(x,y)\right|^2\rmd x\rmd y,\qquad \rho(x,y)=\rme^{\sigma_{\underline{\eta}}(x,y)},\quad \sigma_{\underline{\eta}}(x,y)=\sigma_{\eta_\mathrm{l},\eta_\mathrm{r}}(x)+\sigma_{\eta_\mathrm{b},\eta_\mathrm{t}}(y),
\]
and the associated spaces $H^k_{\underline{\eta}}(\R^2)$.
Clearly, $\displaystyle{\|u\|_{H^k_{\underline{\eta}}(\R^2)}}$ is equivalent to $\displaystyle{\|u(\cdot)\rho(\cdot)\|_{H^k(\R^2)}}$. In other words, multiplication by $\rho(\cdot)$ provides and isomorphism $\displaystyle{H^k_{\underline{\eta}}(\R^2)\to H^k(\R^2)}$. As a consequence, an operator $\mathcal{A}$ of the form \eqref{e:a} is Fredholm on $H^2_{\underline{\eta}}(\R^2)$ if and only if $\mathcal{A}_{\underline{\eta}}$, defined through
\begin{equation}\label{e:wc}
\mathcal{A}_{\underline{\eta}} u=\Delta u + (2\nabla\sigma_{\underline{\eta}}+b)\cdot \nabla u + (\Delta \sigma_{\underline{\eta}}+|\nabla\sigma_{\underline{\eta}}|^2+b\cdot\nabla\sigma_{\underline{\eta}}+a)u,
\end{equation}
is Fredholm as an operator on $H^2(\R^2)$. The product structure in $x$- and $y$-weights in the norms shows that the coefficients of $\nabla{A}_{\underline{\eta}}$ satisfy the uniform limit and smoothness assumptions from Corollary \ref{c:1}. More precisely, we find limiting operators that are obtained from the unweighted limiting operators $\mathcal{A}_\mathrm{b,t,l,r}$ by conjugating with the limiting rate functions 
\begin{equation}
\begin{array}{ll}
\mathrm{l}:\quad\sigma_{\eta_\mathrm{b},\eta_\mathrm{t}}(y)+\eta_\mathrm{l}x,&
\qquad \qquad
\mathrm{b}:\quad\sigma_{\eta_\mathrm{l},\eta_\mathrm{r}}(x)+\eta_\mathrm{b}y,\\
\mathrm{r}:\quad\sigma_{\eta_\mathrm{b},\eta_\mathrm{t}}(y)+\eta_\mathrm{r}x,&
\qquad \qquad
\mathrm{t}:\quad\sigma_{\eta_\mathrm{l},\eta_\mathrm{r}}(x)+\eta_\mathrm{t}y.
\end{array}\label{e:asyw}
\end{equation}
One can now investigate invertibility of the asymptotic operators depending on the weights $\eta_j$. By continuity of the Fredholm index, it does not change as long as invertibility of the asymptotic operators is preserved. In fact, one can even show that dimensions of kernel and cokernel are constant in those connected components of $\underline{\eta}\in\R^4$. To make this precise, first define the bounded embeddings
\[
\iota_{\underline{\eta}\mapsto\underline{\eta'}}:H^k_{\underline{\eta}}(\R^2)\to H^k_{\underline{\eta'}}(\R^2), \quad u\to u,\quad \eta_j\geq\eta_j',j=\mathrm{l,r,t,b}.
\]

\begin{Lemma}[Fredholm properties and weights]\label{Fredholm_persistency_lemma} Let $\mathcal{A}$ be as in \eqref{e:a}, closed and densely defined on $L^2_{\underline{\eta}}(\R^2)$. Then the set of $\underline{\eta}\in \R^4$ for which $\mathcal{A}$ is Fredholm is open. Let $\mathcal{N}$ be a connected subset of the Fredholm region. Then the Fredholm index is constant on $\mathcal{N}$ and the kernels are isomorphic with isomorphism $\iota_{\underline{\eta}\mapsto\underline{\eta'}}$ and its adjoint, respectively, for any $\underline{\eta},\underline{\eta'}\in \mathcal{N}$ satisfying $\eta_j\geq\eta_j'$ for all $j$. 
\end{Lemma}
\begin{Proof}
Clearly, the conjugate operators depend in an analytic fashion on the parameters $\eta_j$ such that the Fredholm index is constant on connected components. In order to show the isomorphism properties of kernel and cokernel, it is enough to vary one component, $\eta_\mathrm{t}$, say, by a small amount, and show that the dimension of the kernel is constant, since the embedding $\iota$ gives a one-to-one map from the kernel in the larger space into the smaller space. The analysis in \cite[\S7.1.3]{kato2013perturbation} then shows that the dimension of the kernel is constant. More precisely, we solve $\mathcal{A}(\underline{\eta})u=0$ near $\underline{\eta}=0$, say. For this, we choose projections $P$ and $Q$ on kernel and cokernel of $\mathcal{A}(0)$, respectively, and decompose $u=Pu+(1-P)u=:z_0+z_\mathrm{h}$, obtaining the equivalent system 
\[
Q\mathcal{A}(\underline{\eta})(z_0+z_\mathrm{h})=0,\qquad (1-Q)\mathcal{A}(\underline{\eta})(z_0+z_\mathrm{h})=0.
\]
We can solve the second equation for $z_\mathrm{h}=B(\underline{\eta})z_0$, with $B(\underline{\eta})$ analytic in $\underline{\eta}$ for $\underline{\eta}\sim 0$, and substitute into the first equation, to obtain a finite-dimensional reduced equation, $A(\underline{\eta})u_0=0$. By construction,  $A(0)= 0$, and, because of the natural inclusion, $A(\underline{\eta})=0$ for $\eta_\mathrm{t}<0$, fixing other weights. By analyticity in $\eta_\mathrm{t}$, we conclude that $A(\underline{\eta})=0$ also for $\eta_\mathrm{t}>0$, locally, as claimed. This proves the lemma. 
\end{Proof}

\begin{Remark}\label{r:1d}
Equivalent results can be derived and are well known in a one-dimensional setting, omitting $y$-derivatives and $y$-dependence in our setup. The operator $\mathcal{A}$  on the real line is then Fredholm when the operators $\mathcal{A}_\mathrm{l/r}$ at $x=\pm\infty$ are invertible; see for instance \cite{palmer,Salamon,ssmorse} for computations of the Fredholm index. 
\end{Remark}

\begin{Proposition}\label{Fredholm_properties:range_of_parameters:2D} The operator 
$$\mathscr{L}: H_{\underline{\eta}}^2(\R^2) \to L_{\underline{\eta}}^2(\R^2)$$
is a Fredholm operator in a connected  open set of weights containing $(\eta_\mathrm{l}, \eta_\mathrm{r}) \in \left((-\delta,0)\cup(0,c_x)\right)\times [0,c_x]$, $\eta_\mathrm{b/t}=0$. In particular, $\mathscr{L}$ is Fredholm for $\eta_j\equiv\eta_0$, $j=\mathrm{b,t,l,r}$, $0<\eta_0$ sufficiently small. 
\end{Proposition}
\begin{Proof} We need to show that the farfield operators $\mathscr{L}_j$, $j=\mathrm{b,t,l,r}$, conjugated with the weights as listed in \eqref{e:asyw}, are bounded invertible. Since the $y$-weights are trivial, the calculation is much simplified. 

We start by considering the linearization at the top, 
\[
\mathscr{L}_\mathrm{t}u=\Delta u + c_x u_x +\mu(x)u-3u_\mathrm{t}^2u,
\]
in the space with weight $\rme^{\sigma_{\eta_\mathrm{l},\eta_\mathrm{r}}}(x)$, $\eta_j\in(0,c_x)$. Fourier transform in $y$ shows that it is sufficient to establish that the spectrum of $\mathscr{L}_\mathrm{t}$ has negative real part. We first consider the case $\eta_j=0$, and then invoke Lemma \ref{Fredholm_persistency_lemma}  and Remark \ref{r:1d} to conclude the general case. From \cite{Monteiro_Scheel}, the operator $\partial_{xx}+c_x\partial_x + \mu(x) -3u_\mathrm{t}^2$ is negative definite at $c_x=0$. Continuity in $c_x$ and Fourier transform in $y$ then readily imply that $\mathscr{L}_\mathrm{t}$ is invertible for all $c_x\geq 0$, small\footnote{This holds true also for finite speeds $c_x$ as seen in the proof of Proposition \ref{p:0}}.

%
%
%
%
%

For nonzero weights, we only need to verify that the linearization at the asymptotic states, 
\[
\mathscr{L}_\mathrm{tl}u=\Delta u + c_x u_x -2u,\qquad \mathscr{L}_\mathrm{tr}=\Delta u + c_x u -u,
\]
are invertible in exponentially weighted spaces with weights $\eta_\mathrm{l/r}$, respectively, which follows readily from Fourier transform. The linearization $\mathscr{L}_\mathrm{b}$ can be shown to be invertible in the same manner.

We next turn to the linearization $\mathscr{L}_\mathrm{l}$,
\[
\mathscr{L}_\mathrm{l}u=\Delta u + c_x u +u-3u_\mathrm{l}^2u.
\]
Conjugating with the exponential weight $\eta_\mathrm{l}$ and using Fourier transform in $x$ gives 
\[
\widehat{\mathscr{L}_\mathrm{l}}(\ell,\eta_\mathrm{l})u=\partial_{yy} u +(\rmi\ell-\eta_\mathrm{l})^2 +  c_x (\rmi \ell -\eta_\mathrm{l}) u +u-3u_\mathrm{l}^2u.
\]
Since the spectrum of the self-adjoint operator $\widehat{\mathscr{L}_\mathrm{l}}(0,0)$ is negative except for the simple eigenvalue $\lambda=0$, we find that $\widehat{\mathscr{L}}_\mathrm{l}(\ell,\eta_\mathrm{l})$ is invertible provided that $r(\ell,\eta_\mathrm{l})=(\rmi\ell-\eta_\mathrm{l})^2+c_x(\rmi\ell-\eta_\mathrm{l}\not\in \R^+$ for all $\ell\in\R$. One finds that the imaginary part vanishes only when $\eta=c_x/2$, which however yields $r(\ell,c_x/2)=-\ell^2-c^2/4<0$. On the other hand, vanishing imaginary part gives $\ell=0$, and $r(0,\eta_\mathrm{l})=\eta_\mathrm{l}^2-c_x\eta_\mathrm{l}$. This quantity is negative for $\eta_x\in (0,c_x)$ and not in the spectrum of $\widehat{\mathscr{L}}_\mathrm{l}(0,0)$ for $\eta_\mathrm{l}<0$, sufficiently small. This shows invertibility of $\mathscr{L}_\mathrm{l}$ with weights $0<\eta_\mathrm{l}<c_x$ as needed. 

Invertibility of $\mathscr{L}_\mathrm{r}$ is easily established in a similar fashion. 
 \end{Proof}

\subsection{Exponential asymptotics in the far field}\label{s:3.2}

We establish (A2) for $c_x$ sufficiently small. From \cite{Monteiro_Scheel}, the linearization at $\Theta$ is invertible on $L^2(\R^2)$ for $c_x=0$, hence for $c_x>0$, small, in the subspace of functions odd in $y$. The results from the previous section can easily be adapted to show that invertibility therefore holds in spaces of functions odd in $y$ with small weights $\eta_j$, $j=\mathrm{l,r,b,t}$, $\eta_\mathrm{b}=\eta_\mathrm{t}$. 

We will now show that the residual $w_0=\Theta(\cdot)-u_\mathrm{ff}(\cdot;0,0)$ is exponentially localized, that is, it belongs to $H^2_\eta(\R^2)$. 
First, consider $w=\Theta_x\chi_{-1}$, where $\chi$ is smooth, $\chi_{-1}=1$ for $x<-2$ and $\chi_{-1}=0$ for $x>-1$. Clearly $w$ is bounded and the residual, $\mathscr{L}w=f$ is odd and an element of $L^2_{\delta,\delta,-\delta',-\delta'}(\R^2)$, for some $\delta,\delta'>0$, since it is bounded and supported in $x\in (-2,-1)$. As a consequence, $w\in H^2_{\delta,\delta,-\delta',-\delta'}(\R^2)$. Integrating,
\[
\Theta(x,y)-u_\mathrm{l}(y;0,0)=\int_{-\infty}^x w(x',y)\rmd x',\qquad x<-2,
\]
gives $\Theta\chi_{-1}\in H^2_{\delta,\delta,\-\delta',-\delta'}(\R^2)$. Inspecting the norms and the support of $\chi_\mathrm{l}$, we immediately see that $(\Theta-u_\mathrm{l}(y;0))\chi_\mathrm{l}\in H^2_\eta(\R^2)$ for some $\eta>0$, small. 

Next, consider $w=\Theta_y$, which solves $\mathscr{L}w=0$. Since $\mathscr{L}$ is Fredholm in on $L^2_{-\delta,-\delta,\delta',\delta'}(\R^2)$ for all $\delta>0,\delta'\sim 0$, and since $w$ is bounded,  $w\in  H^2_{-\delta,-\delta,-\delta',-\delta'}(\R^2)$, for all $\delta,\delta'>0$, we conclude from Lemma \ref{Fredholm_persistency_lemma} that in fact $w\in  H^2_{-\delta,-\delta,\delta',\delta'}(\R^2)$ for all $\delta,\delta'>0$. Integrating, 
\[
\Theta(x,y)-u_\mathrm{t}(x;0)=\int_{\infty}^y w(x,y')\rmd y',
\]
we find that $\Theta\in H^2_{-\delta,-\delta,\delta',\delta'}(\R^2)$, and therefore, inspecting the values of the norm in the sector giving the support of $\chi_\mathrm{t}$, $\chi_\mathrm{t}(\Theta-u_\mathrm{t})\in H^2_\eta(\R^2)$ for $\eta>0$, sufficiently small. 

The estimates for the limits at the bottom and to the right are similar and easier, respectively. This establishes (A2) for $c_x>0$, sufficiently small.

\subsection{The Fredholm index and cokernel}\label{s:3.3}

The linearization $\mathscr{L}$ of \eqref{e:ell} at $\Theta$, 
\[
\mathscr{L}v =\Delta v + c_xv_x   + \mu(x) v - 3\Theta^2v,
\]
clearly leaves invariant subspaces of functions that are even in $y$, or subspaces of functions that are odd in $y$, since coefficients $\Theta(x,y)$ are even in $y$. From \cite[Lem. 4.6]{Monteiro_Scheel}, we conclude that $\mathscr{L}$ is bounded invertible on the odd subspace, $L^2_{\eta,\mathrm{odd}}(\R^2)$, for all $\eta$ with $|\eta|\sim0$ sufficiently small. Next, recall that Fredholm properties are additive for direct sums. In particular, Fredholm indices and dimensions of kernels and cokernels for $\mathscr{L}$ on $L^2_{\eta}(\R^2)$ are the sums of those for $L^2_{\eta,\mathrm{odd}}(\R^2)$ and $L^2_{\eta,\mathrm{even}}(\R^2)$. We therefore only need to show that the kernel of $\mathscr{L}$ is trivial and the cokernel is one-dimensional on $L^2_{\eta,\mathrm{even}}(\R^2)$, $\eta>0$ sufficiently small.

\begin{Lemma}\label{Lemma:almost_trivial_kernel} Consider $h\in  H^2_\eta(\R^2)$, $\eta>0$, even in $y$, $\mathscr{L}h=0$. Then $h=0$. 
\end{Lemma}
\begin{Proof} From Proposition \ref{Fredholm_properties:range_of_parameters:2D} and Lemma  \ref{Fredholm_persistency_lemma}, we conclude that $\displaystyle{e^{c_x x/2}h(x,y) \in L^2(\R^2)}$.  The following  proof mimics an argument found in \cite[\S 4]{berestycki1997further}, also explored in \cite{Kowalczyk}.  Clearly, $\Theta_y$ belongs to the kernel and is even, although not exponentially localized. Since $\Theta_y>0$ we may define 
\[
H = e^{\frac{c_x}{2}x}h, \quad W = e^{\frac{c_x}{2}x}\Theta_y = \partial_y( e^{\frac{c_x}{2}x}\Theta), \quad \sigma =\frac{H}{W} = \frac{h}{\Theta_y}. 
\]
Both $W$ and  $H = \sigma W$ solve $\tilde{\mathscr{L}}v=0$, where
$\tilde{\mathscr{L}}v := \left[\Delta + \left(\mu(x) - 3\Theta^2 - \frac{c_x^2}{4}\right)\right]v.$ We conclude that 
\begin{align}
\sigma W\tilde{\mathscr{L}}(\sigma W) &= \sigma W\Delta(\sigma W) + \sigma W\left(\mu(x) - 3\Theta^2 - \frac{c_x^2}{4}\right)(\sigma W) \nonumber\\
& = \sigma W\left(W \Delta \sigma + 2 \nabla \sigma \nabla W  \right) + \sigma^2 W \left\{\Delta W + \left(\mu(x) - 3\Theta^2 - \frac{c_x^2}{4}\right)W\right\}\nonumber \\
& = \sigma \mathrm{div}\,(W^2 \nabla\sigma).\label{e:div}
\end{align}
Now, define 
\[
\zeta_R = \zeta\left(\frac{\sqrt{x^2 + y^2}}{R}\right),\qquad   \zeta(z) = \left\{\begin{array}{cc}
      1, & \mbox{for} \quad z \leq 1;\\
      0, & \mbox{for} \quad z\geq 2.
      \end{array} \right.
\]
Multiply \eqref{e:div} by $\zeta_R^2$ and integrate to obtain, using the notation $r(x,y)=\sqrt{x^2+y^2}$,
\begin{align*}
\int_{\mathbb{R}^2} \zeta_R^2|\nabla \sigma|^2 W^2 dxdy &= -\frac{2}{R}\int_{\mathbb{R}^2}\zeta_R W^2 \sigma\zeta_R'(\nabla r\cdot \nabla \sigma)dx dy  \nonumber \\
&\leq \frac{2}{R}\left( \int_{R \leq r \leq 2R} \zeta_R^2|\nabla \sigma|^2 W^2 dxdy\right)^{\frac{1}{2}}\left(\int_{R \leq r \leq 2R}|\zeta_R'|^2 |\nabla|r|||\sigma W|^2 \right)^{\frac{1}{2}} \nonumber \\
&\leq \frac{2}{R}\left( \int_{R \leq r \leq 2R} \zeta_R^2|\nabla \sigma|^2 W^2 dxdy\right)^{\frac{1}{2}}\left(\int_{R \leq r \leq 2R}|\zeta_R'|^2 |\nabla|r|||H|^2 \right)^{\frac{1}{2}}.
\end{align*}
Since $|\nabla r| \leq 1$ and  $H\in L^2(\R^2)$, we can choose $R >1$ and obtain the existence of a constant $C$ independent of $R$ such that
\begin{align*}
\int_{\mathbb{R}^2} \zeta_R^2|\nabla \sigma|^2 W^2 dxdy &\leq C^{\frac{1}{2}}\left( \int_{R \leq r \leq 2R} \zeta_R^2|\nabla \sigma|^2 W^2 dxdy\right)^{\frac{1}{2}}.
\end{align*}
from this, we conclude  that  $\int_{\mathbb{R}^2} \zeta_R^2|\nabla \sigma|^2 W^2 dxdy$ is bounded by $C$. Next, letting $R \to \infty$ on the right-hand side we conclude from  Lebesgue's Dominated Convergence Theorem that $ \int_{\mathbb{R}^2} |\nabla \sigma|^2 W^2 dxdy =0.$ Now, since  $W >0$,  we find $\nabla \sigma =0$, which proves that $h=\rho \Theta_y$ for some constant $\rho$. Since $\Theta_y\not\in H^2_\eta(\R^2)$ for $\eta>0$, $\rho=0$ which concludes the proof. 
\end{Proof}
It remains to show that the cokernel of $\mathscr{L}$ is one-dimensional in $L^2_\eta(\R^2)$. We therefore consider the $L^2$-adjoint operator 
\[
\mathscr{L}^* = \Delta - c_x \partial_x  + \mu(x)  - 3\Theta^2,
\]
with domain $H^2_{-\eta}(\R^2)$.  One readily verifies that  $e^{c_x x} \Theta_y \in \ker(\mathscr{L}^*).$ Note that $\mathscr{L}$ and $\mathscr{L}^*$ are conjugate through the multiplication operator $\rme^{c_x x}$. In particular, if $h^* \in \ker\left(\mathscr{L}^*\right)$ then $e^{c_x x}h^* \in \ker\left(\mathscr{L}\right).$  We can therefore use a slight variation of the arguments in Lemma \ref{Lemma:almost_trivial_kernel} to show that the kernel of $\mathscr{L}^*$ is one-dimensional. 
\begin{Lemma}\label{lemma_L_star_kernel}
Let $h^* \in H^2_{-\eta}(\R^2)$, even in $y$,  belong to the kernel of $\mathscr{L}^*$.  Then $h^*$ is a scalar multiple of $\Theta_y$.
\end{Lemma}
\begin{Proof}
Define $H^* = e^{-c_x x/2}h^*$. We have that 
\begin{align*} 
0 = \mathscr{L}^* h^*  = \mathscr{L}^*(e^{c_x x/2}H^*) = e^{c_x\,x/2}\left[\Delta_{x,y} + \left(\mu(x) - 3\Theta^2 - \frac{c_x^2}{4}\right)\right]H^*.
\end{align*}
Now, defining $W^* = e^{c_x x/2}\Theta_y$, we also have that   $\mathscr{L}^*[e^{c_x x/2}W^*]=0$. Next, with 
\[
\sigma^* = \frac{H^*}{W^*} = \frac{e^{-c_xx/2}h^*}{e^{c_x x/2}\Theta_y},
\]
we can follow the proof of  Lemma \ref{Lemma:almost_trivial_kernel} to conclude that $\sigma^*$ is constant a.e, hence 
$h^*=\rho e^{c_x x}\Theta_y$ for some scalar $\rho>0$. 
\end{Proof}
Summarizing, we have established (A3) for $c_x$ sufficiently small. 
\begin{Proposition}[\,(A3) holds for $c_x>0$, small]\label{Theorem:index_linearized_operador}
The operator $\mathscr{L}: H^2_\eta(\R^2) \to L^2_\eta(\R^2)$ is a Fredholm operator with index -1, with trivial kernel, and with cokernel spanned by $\rme^{c_x x} \Theta_y$. 
\end{Proposition}

\begin{Remark}[Spectral flow] One would in general compute the Fredholm index using a spectral flow argument; see for instance \cite{Salamon,ssmorse}.
\end{Remark}

\section{Applications and Discussion}\label{s:5}
We first give brief examples in which we compute the sign of $\frac{\rmd \varphi}{\rmd\alpha}$ from \eqref{e:dp}, Section \ref{s:5.1}. We then discuss our results and possible extensions, also pointing to related results in the literature. 

\subsection{Examples of contact angle selection}\label{s:5.1}
Recall from \eqref{e:dp},\eqref{e:mel}, and \eqref{e:ma}, that
\begin{align*}\label{e:exp}
\frac{\rmd\varphi}{\rmd\alpha}&=-\frac{M_\alpha}{M_\psi},\\
M_\psi&=-c_x\int_{\R^2}\left(\Theta_y\right)^2\rme^{c_x x}\rmd x \rmd y<0,\\
M_\alpha&=-c_\mathrm{n}'(0)M_\psi-\int_\R \rme^{c_x x} \left[G(x,u_\mathrm{t})-G(x,u_\mathrm{b})\right] \rmd x,
\end{align*}
with $G_j'(u)=-g_j(u)$, $j=\mathrm{l/r}$, and $c_\mathrm{n}'(0)$ from \eqref{e:c'}. 

First, consider $g_\mathrm{l}(u)=0$, such that $c_\mathrm{n}(\alpha)\equiv 0$, that is, interfaces in the left half plane $x<0$ do not propagate. Note however that $c_y\neq 0$ in general, when $\psi\neq 0$. For $g_\mathrm{r}(u)=1$, we find $G(x,u)=-u$ and $M_\alpha>0$ in \eqref{e:ma}. Using $M_\psi<0$ from \eqref{e:mel}, we find $\frac{\rmd \varphi}{\rmd\alpha}>0$.

Intuitively, a contact angle greater than $\pi/2$ implies that at the contact line, $x=0$, the interface is propagating downwards, hence at the contact line, the region where $u>0$ is expanding. This aligns well with the intuition where a positive equilibrium state in $x>0$ would facilitate the selection of $u=1$ rather than $u=-1$. 

In this light, it is worth noticing that nonzero contact angles are not created by an imbalance in the energy. in fact, we can choose $g_\mathrm{l}(u)=\frac{1}{2}u^2$, thus retaining the equilibrium state $u_\mathrm{r}(\alpha)\equiv 0$, for all $\alpha$. Yet, 
\[
M_\mathrm{\alpha}=-2\int_\R\rme^{c_x x}u_\mathrm{t}^3\rmd x<0,
\]
since $u_\mathrm{t}=-u_\mathrm{b}$. 

Next, starting with the selection of a contact angle in $x>0$, one can now add relatively small effects in $x<0$, such as $g_\mathrm{l}(u)=\epsilon$, thus changing the speed $c_\mathrm{n}'(0)=\frac{3}{\sqrt{2}}\varepsilon$. 

In fact, these considerations give an interpretation to the two contributions in $M_\alpha$. The second term gives the speed of the contact point between interface and contact line, that is the vertical speed of the point $u=0$ along the quenching line $x=0$. The first term is a simple geometric adjustment to the contact angle such that normal speed in the wake and horizontal propagation of the interface combined result in precisely this vertical propagation.  The contribution to the motion of the contact point, through the integrals of $G$, is indeed exponentially localized near the contact line: the exponential prefactor localizes the effect in $x>-M$, say, and the exponential decay of $u_\mathrm{t/b}$ for $x\to+\infty$ enforces localization in $x<M$, say. 

Slightly generalizing our results, we could have considered $g(x,u)$ converging to $g_\mathrm{l/r}(u)$, exponentially. Choosing $g(x,u)=0$ in $-M<x<0$, $M\gg 1$, $g(x,u)=\varepsilon$ in $x<-M-1$, we see that the contribution of to the integral in the definition of $M_\alpha$ is exponentially small in $M$. With this choice of $g$ in $x<0$ and $g(x,u)=1$ in $x>0$ we can therefore control normal speed $c_\mathrm{n}$ and contact angle $\varphi$ independently, choosing $\varepsilon$ not necessarily small.

\subsection{Summary and future directions}\label{s:5.2}

We presented perturbative results that characterize the creation of interfaces at an internal discontinuity, where system parameters change. At the heart of the analysis is a Fredholm theory that, through a negative Fredholm index, exhibits the necessity of adjusting a farfield matching parameter, naturally chosen as the angle of the interface. The Fredholm analysis and the partition of unity constructions are reminiscent of and to some extend inspired by work on multiple-end solutions in the Allen-Cahn equation; see for instance \cite{pino}. The moving quenching line and the possibility of propagating fronts create however technical differences, such as non-selfadjoint operators. An alternative approach would have adapted the spatial dynamics techniques from \cite{Haragus_Scheel_corner} to this situation, giving of course equivalent results. 

The most natural extension would be to non-small perturbations, preserving the asymptotic monostable and bistable character of the equation, respectively. Results in this global spirit have been obtained in the context of front propagation in homogeneous media, where propagation is accelerated along lines with fast diffusion; see for instance \cite{pauthier}.

Phenomenologically, one can envision more complicated dynamics in the wake. Beyond planar fronts, simple structures known to govern interfacial dynamics are for instance conical fronts \cite{hamel,Haragus_Scheel_corner,Haragus_Scheel_review}, or, in our language here, corners between interfacial lines of different angles. Depending on their horizontal speed of propagation, such corners may or may not interact strongly with the contact line. 

Within the perturbative setup considered here, we would only look at obtuse corners, which propagate at small speeds, hence would not be able to form bound states with the contact line. In order to study such a possibility in more detail, one would therefore want to study small speeds $c_x\gtrsim 0$. In that setting, one would envision the possibility of weakly absorbing contact lines as the dominant structure, similar to the ``holes'' in interfaces constructed in \cite{Haragus_Scheel_corner} or the contact defect structures in \cite{kollar}. In further analogy to \cite{kollar}, see also \cite{ssdefect}, solution constructed here are ``sources'' generating interface, with pointwise transport away from the contact line. The fact that such transport leads to negative Fredholm indices had been noticed in \cite{ssdefect}; see also \cite{ssmorse}.

More basically, in the case of small speeds, the question of interface being ``generated'' at the boundary becomes more subtle, since interface propagation at angles with large enough speed  may effectively lead to interface being absorbed in the boundary. 

Comparing with the results in \cite{Monteiro_Scheel}, one would wish to extend the results here to situations periodic in $y$, or to more general equations such as the Cahn-Hilliard equation. Results on such periodic configurations, with two-dimensional structure have recently been obtained in \cite{gohscheel} for the Swift-Hohenberg equation, again in a perturbative setting.

%
%
\end{document}